\pdfoutput=1

\documentclass[letterpaper, 11pt]{article}
\usepackage{amssymb,amsmath}
\usepackage{epsf}
\usepackage{times,bm,mathrsfs}  % Fonts
\usepackage{amsmath}
\usepackage{amssymb}
\usepackage{amsfonts}
\usepackage{mathrsfs}
\usepackage{bbm}
\usepackage{color}
\usepackage{graphicx}
\usepackage{amscd}
\usepackage{epsfig}

\newenvironment{proof}{\noindent{\textbf{Proof:}}}{$\blacksquare$\vskip\belowdisplayskip}
\newcommand{\mcfn}{{\cal M}^{\mathrm{CFN}}_{f,g}}
\newcommand{\anc}[2]{\widehat{\mathrm{Anc}}_{#1,#2}}
\newcommand{\ancname}{\widehat{\mathrm{Anc}}}
\newcommand{\Maj}{\widehat{\mathrm{Maj}}}
\newcommand{\MajCorr}{{\mathrm{MajCorr}}}
\newcommand{\hatd}{\widehat{\mathrm{Int}}}
\newcommand{\internal}[4]{\hatd\left(#1,#2; #3, #4\right)}
\newcommand{\logcorr}{\widehat{\mathrm{Dist}}}
\newcommand{\distance}[2]{\logcorr\left(#1, #2\right)}
\newcommand{\true}{\textsc{true}}
\newcommand{\false}{\textsc{false}}
\newcommand{\distest}{\textsc{DistanceEstimate}}
\newcommand{\recseq}[1]{\widehat{\sigma}_{#1}}
\newcommand{\itemname}[1]{$\mathrm{[#1]}$}
\newcommand{\onetwo}{\phi} % don't use i !
\newcommand{\issplit}{\textsc{IsSplit}}
\newcommand{\isshort}{\textsc{IsShort}}
\newcommand{\iscollision}{\textsc{IsCollision}}
\newcommand{\bubble}{\textsc{RemoveCollision}}

\newcommand{\cherryid}{\textsc{LocalCherry}}
\newcommand{\fakecherry}{\textsc{DetectCollision}}

\newcommand{\quartet}{\qcal}
\newcommand{\hatquartet}{\widehat{\quartet}}
\newcommand{\forest}[1]{\mathcal{F}_{#1}}
\newcommand{\roots}[1]{\rho(\mathcal{F}_{#1})}

\newcommand{\metric}[1]{\mathcal{D}_{#1}}
\newcommand{\maxmetric}[1]{\overline{\metric{}}_{#1}}

\newcommand{\ischerry}{\mathrm{IsCherry}}
\newcommand{\neighbors}{\mathcal{N}}
\newcommand{\parent}[1]{\mathrm{Parent}_{#1}}
\newcommand{\sister}[1]{\mathrm{Sister}_{#1}}
\newcommand{\forestno}{\mathcal{F}}
\newcommand{\metricno}{\mathcal{D}}
\newcommand{\rootsno}{\rho(\forestno)}
\newcommand{\hascollision}{\mathrm{HasCollision}}
\newcommand{\erada}{\widehat{R}_{\mathrm{acc}}}
\newcommand{\eradc}{\widehat{R}_{\mathrm{col}}}

\newcommand{\bias}{B}
\newcommand{\tol}{\eps}
\newcommand{\distmet}{\textsc{DistortedMetric}}
\newcommand{\remforestno}{\widetilde{\forestno}}
\newcommand{\remforest}[1]{\widetilde{\mathcal{F}}_{#1}}

\definecolor{Red}{rgb}{1,0,0}
\definecolor{Blue}{rgb}{0,0,1}
\definecolor{Olive}{rgb}{0.41,0.55,0.13}
\definecolor{Green}{rgb}{0,1,0}
\definecolor{MGreen}{rgb}{0,0.8,0}
\definecolor{DGreen}{rgb}{0,0.55,0}
\definecolor{Yellow}{rgb}{1,1,0}
\definecolor{Cyan}{rgb}{0,1,1}
\definecolor{Magenta}{rgb}{1,0,1}
\definecolor{Orange}{rgb}{1,.5,0}
\definecolor{Violet}{rgb}{.5,0,.5}
\definecolor{Purple}{rgb}{.75,0,.25}
\definecolor{Brown}{rgb}{.75,.5,.25}
\definecolor{Grey}{rgb}{.5,.5,.5}
\definecolor{Black}{rgb}{0,0,0}

\newcommand{\ecal}{\mathcal{E}}
\newcommand{\fcal}{\mathcal{F}}

\newcommand{\lcal}{\mathcal{L}}
\newcommand{\ncal}{\mathcal{N}}

\newcommand{\qcal}{\mathcal{Q}}

\newcommand{\vcal}{\mathcal{V}}

\newcommand{\ycal}{\mathcal{Y}}

\newcommand{\real}{\mathbb{R}}

\newcommand{\eps}{\varepsilon}

\newcommand{\ind}{\mathbbm{1}}

\newcommand{\bdm}{\begin{displaymath}}
\newcommand{\edm}{\end{displaymath}}
\newcommand{\bea}{\begin{eqnarray*}}
\newcommand{\eea}{\end{eqnarray*}}
\newcommand{\bean}{\begin{eqnarray}}
\newcommand{\eean}{\end{eqnarray}}

\newcommand{\prob}{\mathbb{P}}
\newcommand{\expec}{\mathbb{E}}

\newcommand{\poly}{\mathrm{poly}}

\newcommand{\quantum}{\Delta}

\newcommand{\forestz}[1]{(\mathcal{F}_0)_{#1}}
\newcommand{\rootsz}[1]{\rho(\forestz{#1})}
\newcommand{\remforestz}[1]{(\widetilde{\mathcal{F}}_0)_{#1}}

\newtheorem{theorem}{Theorem}
\newtheorem{corollary}[theorem]{Corollary}

\newtheorem{proposition}{Proposition}[section]
\newtheorem{definition}[proposition]{Definition}

\newtheorem{lemma}[proposition]{Lemma}

\numberwithin{figure}{section}

\newcommand{\ignore}[1]{}

\newcommand{\sign}{{\mbox{{\rm sign}}}}

\renewcommand{\poly}{{\mbox{{\rm poly}}}}

 \def\eps{\varepsilon} \def\E{{\bf{E}}}

\def\R{\hbox{I\kern-.2em\hbox{R}}}

\def\|{\, | \, }
\def\v0{{\bf 0}}

\def\0{\hat{0}}
\def\1{\hat{1}}

\def\path{{\tt path}}

\def\phi{\varphi}

\def\be{\begin{equation}}
\def\ee{\end{equation}}

\def\L{{\cal L}}

\def\part{{\cal P}}
\def\vert{{\cal V}}
\def\edges{{\cal E}}

\def\R{\mathcal{R}}

\def\L{\mathcal{L}}

\def\eps{\varepsilon}

\author{
Constantinos Daskalakis\thanks{CSAIL,
MIT. Work done at U.C Berkeley and Microsoft Research. Partially supported by CIPRES (NSF ITR grant \# NSF EF 03-31494). Email: {\tt costis@csail.mit.edu}.}
\and
Elchanan Mossel\thanks{Statistics and Computer Science,
U.C. Berkeley, and
Weizmann Institute of Science, Rehovot, Israel.
Supported by a
Miller fellowship in Statistics and Computer Science, by a Sloan
fellowship in Mathematics and by NSF grants DMS-0504245,
DMS-0528488 DMS-0548249 (CAREER), ONR grant N0014-07-1-05-06 and ISF grant 1300/08. Email: {\tt mossel@stat.berkeley.edu}.}
\and S\'ebastien Roch\thanks{Department of Mathematics, UCLA. Work at U.C. Berkeley and Microsoft Research. Partially supported by CIPRES (NSF ITR grant \# NSF EF 03-31494). Email: {\tt sebastien.roch@gmail.com}.}}
\title{Evolutionary Trees and the Ising Model on the
Bethe Lattice:\\ A Proof of Steel's Conjecture}

\begin{document}

\maketitle

\thispagestyle{empty}

\begin{abstract}
A major task of evolutionary biology is the
reconstruction of phylogenetic trees from molecular data. The
evolutionary model is given by a Markov chain on
a tree. Given samples from the leaves of the Markov chain,
the goal is to reconstruct the
leaf-labelled tree.

It is well known that in order to
reconstruct a tree on $n$ leaves,
sample sequences
of length $\Omega(\log
n)$ are needed. It was conjectured by M.~Steel that for the CFN/Ising
evolutionary model, if the mutation probability on all edges of
the tree is less than $p^{\ast} = (\sqrt{2}-1)/2^{3/2}$,
then the
tree can be recovered from sequences of length $O(\log n)$. The value
$p^{\ast}$ is given by the transition point for the extremality of the
 free Gibbs measure for the Ising model on the binary tree. Steel's
conjecture was proven by the second author in the special case where the tree
is ``balanced.'' The second author also proved that if all edges
have mutation probability larger than $p^{\ast}$ then the length
needed is $n^{\Omega(1)}$.

Here we show that Steel's conjecture holds true for general trees by giving a
reconstruction algorithm that recovers the tree from $O(\log n)$-length sequences
when the mutation probabilities are discretized and less than $p^\ast$.
Our proof and results demonstrate that extremality of the free
Gibbs measure on the infinite binary tree, which has been studied
before in probability, statistical physics and computer science,
determines how distinguishable are Gibbs measures on finite binary
trees.
\end{abstract}

%\vspace{1mm}
%\noindent
%{\bf Categories and Subject Descriptors:}
%F.2 {[Theory of Computation]}: {Analysis of Algorithms and Problem Complexity};
%J.3 {[Computer Applications]}: {Life and Medical Sciences}---{\em Biology and genetics}.

%\vspace{1mm} \noindent {\bf General Terms:} Phylogenetics,
%Efficient Algorithms.
% EM: Changed general terms

%\vspace{1mm}
%\noindent

\section{Introduction}
In this paper we prove a central conjecture in mathematical phylogenetics~\cite{Steel:01}:
%{\red In this paper we make significant progress
%towards proving a central conjecture in computational
%phylogenetics~\cite{Steel:01}}:
We show that every phylogenetic tree
with short, discretized edges on $n$ leaves can be reconstructed from sequences
of length $O(\log n)$, where by short we mean that the mutation probability
on every edge is bounded above by the critical transition probability for the extremality
of the Ising model on the infinite binary tree.

This result establishes that the extremality of
the free Gibbs measure for the Ising model on the infinite binary
tree, studied in probability, statistical physics and computer
science, determines the sampling complexity of the phylogenetic
problem, a central problem in evolutionary biology. We proceed with
background on the phylogenetic problem, on the
reconstruction problem and a statement of our results.\\

%{\green [SR: The Phylogeny Background should be updated to be more ``abstract'' and
%probabilist-friendly. My suggestion is to remove any mention of DNA, character, Purines, Pyrimidies, A,C,G,T, mutation, etc.
%Maybe just talk about inferring Ising models on trees from samples at the leaves.
%I would also remove completely the Potts case.
%We can come back at the end of the introduction
%to a short section explaining the motivation for these definitions.]}

%{\blue [EM: Agree. I commented most of this stuff out.]}

\noindent{\bf Phylogenetic Background.}
Phylogenies are used in evolutionary biology to model the stochastic
evolution of genetic data on the ancestral tree relating a group of species.
The leaves of the tree correspond to (known)
extant species. Internal nodes represent
extinct species. In particular the root of the tree represents the most
recent ancestor of all species in the tree. Following
paths from the root to the leaves, each bifurcation indicates
a speciation event whereby two new species are created
from a parent. We refer the reader
to~\cite{SempleSteel:03}
for an excellent introduction to phylogenetics.
The underlying assumption is that genetic information evolves from
the root to the leaves according to a Markov model on the tree.
It is further assumed that this process is repeated independently a number of times
denoted $k$.
Thus each node of the tree is associated with a sequence of length $k$.
The vector of the $i$'th letter of all sequences at the leaves is called the
$i$'th {\em character}. One of the major tasks in
molecular biology, the {\em reconstruction of phylogenetic trees},
is to infer the topology of the tree from the characters
at the leaves.

%This genetic information may consist of DNA sequences, proteins,
%etc. Suppose for example that the genetic data consists of
%(aligned) DNA sequences and let us follow the evolution of the
%first letter in all sequences. This collection, named the first
%{\em character}, evolves according to Markov transition matrices
%on the edges. The root is assigned one of the four letters $A,C,G$
%and $T$. Then this letter evolves from parents to descendants
%according to the Markov matrices on the edges connecting them.

%The vector of the $i$'th letter of all sequences is called the
%$i$'th {\em character}. It is further assumed that the characters
%are i.i.d. random variables. In other words, each site in a DNA
%sequence is assumed to mutate independently from its neighbors
%according to the same mutation mechanism. Naturally, this is an
%over-simplification of the underlying biology. Nonetheless, the
%model above may be a good model for the evolution of some DNA
%subsequences and is the most popular evolution model in molecular
%biology, see e.g.~\cite{Felsenstein:04,SempleSteel:03}.
%One of the major tasks in
%molecular biology, the {\em reconstruction of phylogenetic trees},
%is to infer the topology of the (unknown) tree from the characters
%(sequences) at the leaves (extant species).

In this paper we will be mostly interested in two evolutionary models,
the so-called Cavender-Farris-Neyman (CFN)~\cite{Cavender:78,Farris:73,Neyman:71} and Jukes-Cantor
(JC)~\cite{JukesCantor:69} models.
In the CFN model the states at the nodes of the tree
are $0$ and $1$ and their a
priori probability at the root is uniform. To
each edge $e$ corresponds a {\em mutation probability} $p(e)$ which is the
probability that the state changes along the edge $e$.
Note that this model is identical to the free Gibbs measure of the Ising model on the tree.
See~\cite{Lyons:89}.
In the JC model the states are $A$, $C$, $G$ and $T$ with a priori
probability $1/4$ each. To each edge $e$ corresponds a mutation
probability $p(e)$ and it is assumed that every state transitions with
probability $p(e)$ to each of the other states.
This model is equivalent to the ferromagnetic Potts model on the tree.

\noindent{\bf Extremality and the Reconstruction Problem.} A
problem that is closely related to the {\em phylogenetic problem}
is that of inferring the {\em ancestral state}, that is, the
state at the root of the tree, given the states at the
leaves. This problem was studied earlier in statistical physics,
probability and computer science under the name of {\em
reconstruction problem}, or {\em extremality of the free Gibbs
measure}. See~\cite{Spitzer:75,Higuchi:77,Georgii:88}. The
reconstruction problem for the CFN model was analyzed
in~\cite{BlRuZa:95,EvKePeSc:00,Ioffe:96a,BeKeMoPe:05,MaSiWe:04,BoChMoRo:06}.
In particular, the role of the reconstruction problem in the
analysis of the mixing time of Glauber dynamics on trees was
established in~\cite{BeKeMoPe:05,MaSiWe:04}.

Roughly speaking, the reconstruction problem is {\em solvable}
when the correlation between the root and the leaves persists no
matter how large the tree is. When it is unsolvable, the
correlation decays to $0$ for large trees. The results
of~\cite{BlRuZa:95,EvKePeSc:00,Ioffe:96a,BeKeMoPe:05,MaSiWe:04,BoChMoRo:06}
show that for the CFN model, if for all $e$ it holds that $p(e)
\leq p_{\max} < p^{\ast}$, then the reconstruction problem is
solvable, where
\[
p^{\ast} = \frac{\sqrt{2}-1}{\sqrt{8}} \approx 15\%.
\]
If, on the other hand, for all $e$ it holds that
$p(e) \geq p_{\min} > p^{\ast}$ and the tree is balanced in the
sense that all leaves are at the same distance from the root, then the
reconstruction problem is {\em unsolvable}. Moreover in this case, the
correlation between the root state and any function of the character
states at the leaves decays as $n^{-\Omega(1)}$.

\noindent{\bf Our Results.}
M.~Steel~\cite{Steel:01} conjectured that when
$0 < p_{\min} \leq p(e) \leq p_{\max} < p^{\ast}$
for all edges $e$, one can reconstruct with
high probability the phylogenetic tree from $O(\log n)$ characters.
Steel's insightful conjecture suggests that there are deep connections
between the reconstruction problem and phylogenetic reconstruction.

This conjecture has been proven to hold for trees where all the
leaves are at the same graph distance from the root---the so-called
``balanced'' case---in~\cite{Mossel:04a}. It is also shown there that the number of
characters needed when $p(e) \geq p_{\min} > p^{\ast}$, for all $e$,
is $n^{\Omega(1)}$. The second result intuitively follows from the
fact that the topology of the part of the tree that is close to
the root is essentially independent of the characters at the
leaves if the number of characters is not at least
$n^{\Omega(1)}$.

The basic intuition behind Steel's conjecture is that: since in
the regime where $p(e) \le p_{\max} < p^{\ast}$ there is no decay
of the quality of reconstructed sequences, it should be as easy to
reconstruct deep trees as it is to reconstruct shallow trees.
In~\cite{ErStSzWa:99a} (see also~\cite{Mossel:07}) it is shown
that ``shallow'' trees can be reconstructed from $O(\log n)$
characters if all mutation probabilities are bounded away from $0$
and $1/2$ (the results of~\cite{ErStSzWa:99a} also show that in
this regime, any tree can be recovered from sequences of
polynomial length). The same high-level reasoning has also yielded
a complete proof that $O(\log n)$ characters suffice for a
percolation-type mutation model when all edges are
short~\cite{MosselSteel:04a}. See~\cite{MosselSteel:04a} for details.

%{\red Here we make significant progress towards
%a complete proof of Steel's conjecture}.
Here we prove Steel's conjecture for general trees under the assumption that the mutation probabilities 
are discretized.
We show that,
if $0 < p_{\min} \leq p(e) \leq p_{\max} < p^{\ast}$ for all edges
$e$ of the tree, then the tree can be reconstructed from
$c(p_{\min},p_{\max})(\log n + \log{1/\delta})$ characters with
error probability at most $\delta$.
The discretization assumption amounts to assuming that all edge lengths are a multiple of a small constant $\Delta$.
(See below for a formal definition of ``edge length.'')
%We prove our result under the assumption that the $p(e)$'s
%are ``discretized.'' [SR: We may want to say more precisely what we assume but that
%would imply introducing distances.]
This result further implies that sequences of logarithmic length suffice
to reconstruct phylogenetic trees in the Jukes-Cantor model, when
all the edges are sufficiently short.

%{\bf\red [SR: Rephrased paragraph below.]}
Compared to~\cite{Mossel:04a}, our main technical contribution is the design and analysis of
a tree-building procedure that uses only ``local information'' to build parts of the tree,
while also maintaining ``disjointness'' between the different reconstructed subtrees.
The disjointness property is crucial in order to maintain the conditional independence
properties of the Gibbs distribution on the original tree---for the purpose of performing estimation
of the states at internal nodes of the tree.
Note that in the balanced case of~\cite{Mossel:04a} this property can be achieved in a straightforward
manner by building the tree ``level-by-level.''

%{\green [SR: Add a paragraph about our techniques and how they differ from
%[Mossel'04]. Either here or somewhere else in the introduction. It may actually
%be a good idea to have an actual subsection describing our technical contributions.]}

\subsection{Formal Definitions and Main Results} \label{sec:definitions}

\paragraph{Trees and Path Metrics.} Let $T$ be a tree. Write $\vert(T)$ for the nodes of $T$,
$\edges(T)$ for the edges of $T$ and $\L(T)$ for the leaves of
$T$. If the tree is rooted, then we denote by $\rho(T)$ the root
of $T$. Unless stated otherwise, all trees are assumed to be {\em
binary} (all internal degrees are $3$) and it is further assumed
that $\L(T)$ is labelled.

Let $T$ be a tree equipped with a length function on its edges, $d
: \edges(T) \to \real_{+}$. The function $d$ will also denote the induced path
metric on $\vert(T)$: $d(v,w) = \sum \{d(e) : e \in
\path_T(v,w)\}$, for all $v,w \in \vert(T)$, where $\path_T(x,y)$
is the path (sequence of edges) connecting $x$ to $y$ in $T$.

We will further assume below that the length of all edges
is bounded between $f$ and $g$ for all $e \in \edges(T)$. In other words, for
all $e \in \edges(T)$,
$f \leq d(e) \leq g.$

\paragraph{Markov Model of Evolution.} 
The evolutionary process is determined by a rooted tree $T=(V,E)$
equipped with a path metric $d$ and
a {\em rate matrix} $Q$.
We will be mostly interested in the case where
$Q = \left( \begin{smallmatrix} -1 & 1 \\ 1 & -1 \end{smallmatrix} \right)$
corresponding to the CFN model and in the case where $Q$ is a $4
\times 4$ matrix given by $Q_{i,j} = 1 - 4 \cdot \ind\{i=j\}$
corresponding to the Jukes-Cantor model.
To edge $e$ of length $d(e)$ we associate the transition matrix
$M^e = \exp(d(e) Q)$.

In the evolutionary model on the tree $T$ rooted at $\rho$ each vertex
iteratively chooses its state from the state at its parent by an
application of the Markov transition rule $M^e$, where $e$ is the
edge connecting it to its parent. We assume that all edges in $E$
are directed away from the root. Thus the probability distribution
on the tree is the probability distribution on $\{0,1\}^V$
($\{A,C,G,T\}^V$) given by 
$$ \overline{\mu}[\sigma] =
\pi(\sigma(\rho)) \prod_{(x \to y) \in E} M^{(x\rightarrow
y)}_{\sigma(x),\sigma(y)}, $$ 
where $\pi$ is given by the uniform
distribution at the root, so that $\pi(0) = \pi(1) = 1/2$ for the
CFN model and $\pi(A) = \pi(C) = \pi(G) = \pi(T) = 1/4$ for the JC
model. We let the measure $\mu$ denote the marginal of
$\overline{\mu}$ on the set of leaves which we identify with
$[n]=\{1,\ldots,n\}$. Thus 
$$ \mu(\sigma) = \sum \{ \overline{\mu}(\tau) : \forall
i \in [n], \tau(i) = \sigma(i) \}. $$ 
The measure $\mu$ defines the
probability distribution at the leaves of the tree.

We note that both for the CFN model and for the JC model, the
transition matrices $M^e$ admit a simple alternative 
representation. For the CFN model, with probability $p(e) =
(1-\exp(-2d(e)))/2$, there is a transition and, otherwise, there is
no transition. Similarly for the JC model with probability $p(e) =
(1-\exp(-4d(e)))/4$ each of the three possible transitions occur. In
particular, defining
\begin{equation} \label{eq:def_g}
g^{\ast} = \frac{\ln 2}{4},
\end{equation}
we may
formulate the result on the reconstruction problem for the phase
transition of the CFN model as follows:
``If $ d(e)  \leq g < g^{\ast}$ for all $e$ then the
reconstruction problem is solvable.''

\paragraph{Phylogenetic Reconstruction Problem.} We will be interested in
reconstructing phylogenies in the regime
where the reconstruction problem is solvable.
The objective is to reconstruct the underlying tree
whose internal nodes are unknown from the collection of sequences
at the leaves. Since for both the CFN model and the
JC model, the distribution $\overline{\mu}[\sigma]$ described
above is independent of the location of the root we can only aim
to reconstruct the underlying unrooted topology.
Let ${\cal T}$ represent the set of all {\em binary topologies}
(that is, unrooted undirected binary trees) and $\mcfn$
the family of CFN transition matrices, as described above, which
correspond to distances $d$ satisfying:
\begin{equation*}
0 < f \leq d \leq g < g^{\ast},
\end{equation*}
where $g^{\ast}$ is given by (\ref{eq:def_g}) and $f$ is an
arbitrary positive constant. Let ${\cal T} \otimes \mcfn$
denote the set of all unrooted phylogenies, where
the underlying topology is in ${\cal T}$ and all transition matrices
on the edges are in $\mcfn$. Rooting $T \in {\cal
T} \otimes \mcfn$ at an arbitrary node, let $\mu_T$
be the measure at the leaves of $T$ as described above.
% EM: Changed wording
It is well known, e.g.~\cite{ErStSzWa:99a,Chang:96}, that different
elements in ${\cal T} \otimes \mcfn$ correspond to
different measures; therefore we will identify measures with their
corresponding elements of ${\cal T} \otimes \mcfn$.
We are interested in finding a (efficiently computable) map $\Psi$
such that $\Psi(\sigma^1_{\partial},\ldots,\sigma^k_{\partial})
\in {\cal T}$, where $\sigma_{\partial} = \left(\sigma^i_{\partial}\right)_{i=1}^k$
are $k$ characters at the leaves of the tree. Moreover, we require
that for every distribution $\mu_T \in {\cal T} \otimes \mcfn$, if the characters
$\sigma^1_{\partial},\ldots,\sigma^k_{\partial}$ are generated
independently from $\mu_T$, then with high probability
$\Psi(\sigma^1_{\partial},\ldots,\sigma^k_{\partial}) = T$. The
problem of finding an efficiently computable map $\Psi$
is called the {\em phylogenetic reconstruction
problem} for the CFN model. The phylogenetic reconstruction
problem for the JC model is defined similarly.
In~\cite{ErStSzWa:99a}, it is shown that there exists a polynomial-time
algorithm that reconstructs the topology from $k =
\poly(n, \log{1/\delta})$ characters, with probability of error $\delta$.
Our results are the following.
We first define a subset of ${\cal T} \otimes \mcfn$.
In words, the $\quantum$-Branch Model ($\quantum$-BM) is a subset of ${\cal T} \otimes \mcfn$
where the edge lengths $d(e)$, $e\in E$, are discretized. This
extra assumption is made for technical reasons. See Section~\ref{section:beyond}.
\begin{definition}[$\quantum$-Branch Model]
Let $\quantum > 0$. We denote by $\ycal[\quantum]$ the subset
of ${\cal T} \otimes \mcfn$ where all $d(e)$'s are multiples of $\quantum$.
We call $\ycal[\quantum]$ the $\quantum$-Branch Model ($\quantum$-BM).
\end{definition}
\begin{theorem}[Main Result]\label{thm}
Consider the $\quantum$-Branch Model above with: $0
< f \leq g < g^{\ast}$ and $\quantum > 0$.
Then there exists a polynomial-time
algorithm that reconstructs the topology of the tree from $k
= c(f,g,\quantum)(\log n + \log{1/\delta})$ characters with error
probability at most $\delta$. In particular,
$$c(f,g,\quantum)=\frac{c(g)}{\min\{\quantum^2, f^2\}}.$$
Moreover, the value $g^*$
given by (\ref{eq:def_g}) is tight.
\end{theorem}
%{\bf\red [SR: Added refs in Corollary below.]}
\begin{corollary}[Jukes-Cantor Model] \label{cor}
Consider the JC model
on binary trees where all edges satisfy
\[
0 < f \leq g < g_{\text{\tiny JC}}^{\ast},\text{ where
}g_{\text{\tiny JC}}^{\ast}:=g^{\ast}/2.
\]
Under the $\quantum$-BM, there exists a polynomial-time algorithm that reconstructs
the topology of the tree from $c'(f,g,\quantum)(\log n + \log{1/\delta})$
characters with error probability at most $\delta$ and
$$c'(f,g,\quantum)=\frac{c'(g)}{\min\{\quantum^2,f^2\}}.$$
The value $g_{\text{\tiny JC}}^{\ast}$ corresponds to the so-called
Kesten-Stigum bound~\cite{KestenStigum:66} on the reconstruction threshold, which has been conjectured
to be tight for the JC model~\cite{MezardMontanari:06, Sly:08}.
\label{cor:JukesCantor}
\end{corollary}
%{\blue [EM: Commented out info about percentage]}
%Note that in the Jukes-Cantor case we allow up to
%$\approx 22\%$ nucleotide changes on average \emph{per branch}.
%In particular, we believe our result is of potential practical interest.
Theorem~\ref{thm} and Corollary~\ref{cor} extend also to cases
where the data at the leaves is given with an arbitrary level of
noise. For this \emph{robust phylogenetic reconstruction problem}
both values $g^{\ast}$ and $g_{\text{\tiny JC}}^{\ast}$ are tight.
See~\cite{JansonMossel:04}.

The results stated here were first reported without proof in \cite{DaMoRo:06}.
Note that in~\cite{DaMoRo:06} the result were
%mistakenly
stated without
the discretization assumption which is in fact needed for the final step of the proof.
This is further explained in subsection~\ref{section:beyond}.

%{\green
%\subsection{Biological Motivation}

%[SR: To be written.]}

\subsection{Organization of the Paper}

Roughly speaking, our reconstruction algorithm has two main components.
First, the {\em statistical component} consists in
\begin{enumerate}
\item \itemname{Ancestral\ Reconstruction} the reconstruction of sequences at internal
nodes;
\item \itemname{Distance\ Estimation} the estimation of distances
between the nodes.
\end{enumerate}
The former---detailed in Section~\ref{sec:ancestral}---is borrowed from~\cite{Mossel:04a}
where Steel's conjecture is proved for the
special case of balanced trees.
In general trees, however, distance estimation is complicated by
nontrivial correlations between reconstruction biases.
We deal with these issues in Section~\ref{sec:stattech}.
%On the other hand,

Second, the {\em combinatorial component} of the algorithm---which uses quartet-based
ideas from phylogenetics and is significantly
more involved than~\cite{Mossel:04a}---is detailed in
Sections~\ref{sec:combtech} and~\ref{sec:algvars}. A full description of the algorithm
as well as an example of its execution can also be found in Section~\ref{sec:algvars}.
Proof of the correctness
of the algorithm is provided in Section~\ref{sec:analysis} and~\ref{sec:union}.

%On a first reading, the reader may want to skip the proofs in Sections~\ref{sec:stattech}
%and~\ref{sec:combtech}
%since they are rather straightforward and are not essential to understand the algorithm
%and the proof of its correctness.

%The paper is organized as follows. We introduce basic statistical and
%combinatorial tools in Sections~\ref{sec:stattech} and~\ref{sec:combtech}.
%We then give a complete description of the reconstruction algorithm
%in Section~\ref{sec:overview}. In Section~\ref{sec:analysis}
%we provide an analysis of its correctness.

%\section{Algorithm}\label{sec:overview}

\subsection{Notation}

Throughout we fix
$0 < f \leq g < g' < g^*$, $\quantum > 0$, and $\gamma > 3$. By definition
of $g^*$, we have $2e^{-4g}>1$ and $2e^{-4g'}>1$. We let
$\theta = e^{-2g}$ and $\theta' = e^{-2g'}$.

%\section{Statistical Estimators}\label{sec:stattech}

%We begin with a description of the statistical techniques used in the reconstruction
%algorithm.

\section{Reconstruction of Ancestral Sequences}\label{sec:ancestral}

%{\green [SR: I think it would be a nightmare for us and the readers to remove
%$\theta$ here. We should at least write the expression of $\theta$ as a function of $d$
%instead of $p$ and explain carefully the connection between them to help the reader.
%For example, we should take the time to explain what the various conditions on $\theta$
%mean in terms of $d$.]}

\paragraph{Ancestral Reconstruction.} In this section we state the results of~\cite{Mossel:04a} on ancestral
reconstruction using recursive majority and we briefly
explain how these results are used in the reconstruction
algorithm.
Following~\cite{Mossel:04a}, we use state values
$\pm 1$ instead of $0/1$. Furthermore, we use the parameter
$\theta(e) = 1-2p(e) = e^{-2 d(e)}$.
Note that $\theta(e)$ measures the correlation between the states
at the endpoints of $e$.
Because the CFN model is ferromagnetic, we have $0 \leq \theta(e) \leq 1$.
In terms of $\theta$ we have
reconstruction solvability whenever $\theta(e) \geq \theta >
\theta_{\ast}$, for all edges $e$, where the value $\theta_{\ast}$
satisfies $2 \theta_{\ast}^2 = 1$.

For the CFN model both the majority algorithm~\cite{Higuchi:77}
and the recursive majority algorithm~\cite{Mossel:98}\footnote{See below for a definition of these estimators.} are
effective in reconstructing the root value. Note that for other
models, in general, most simple reconstruction algorithms are not
effective all the way to the reconstruction
threshold~\cite{Mossel:01,MosselPeres:03,JansonMossel:04}.
However, as noted in~\cite{Mossel:04a}, there is an important
difference between the two reconstruction algorithms when
different edges $e$ have different values of $\theta(e)$. Suppose
that $\theta(e)
> \theta'
> \theta^{\ast}$ for all edges $e$.
Then, the recursive majority function is effective
in reconstructing the root value with probability bounded away
from $1/2$ (as a function of $\theta'$). On the other hand, it is easy
to construct examples where the majority function reconstructs the
root with probability tending to $1/2$ as the tree size increases.

The difference between the two algorithms can be roughly stated as
follows. When different edges have different $\theta$-values,
different ``parts'' of the tree should have different weights in
calculating the reconstructed value at the root. Indeed, when all
$\theta$-values are known, an appropriate weighted majority function can
estimate the root value correctly with probability bounded away from
$1/2$~\cite{MosselPeres:03}. However, when the
$\theta$-values are unknown, using  uniform weights may result in an
arbitrarily inaccurate estimator.

Recursive majority, on the other hand, does not require
knowledge of the $\theta$-values to be applied successfully. This
essentially follows from the fact that the majority function is
``noise-reducing'' in the following sense. Suppose $\theta' >
\theta^{\ast}$. Then, as we shall see shortly, there exists an
integer $\ell$ and noise level $q <1/2$ such that majority on the
$\ell$-level binary tree has the following property:
if all leaf values are given with stochastic noise at most
$q$, then the majority of these values differs from the actual root state with probability at most $q$. 
Therefore,
recursive application of the majority function---$\ell$ levels
at a time---results in an estimator whose error is at most $q$ for
{\em any} number of levels.

\paragraph{Properties of Recursive Majority.} We proceed with a formal definition
of recursive majority. Let $\Maj :
\{-1,1\}^d \to \{-1,1\}$ be the function defined as follows
$$\Maj(x_1,\ldots,x_d) =  \sign \left(\sum_{i=1}^d x_i + 0.5 \omega\right),$$
where $\omega$ is $\pm 1$ with probability $1/2$ independently of
the $x_i$'s. In other words, $\Maj$ outputs the majority value of
its input arguments, unless there is a tie in which case it
outputs $\pm 1$ with probability $1/2$ each. For consistency, we denote
all statistical estimators with a hat. Note in particular that
our notation differs from~\cite{Mossel:04a}.

Next we define the ``noise-reduction'' property of majority. The function
$\eta$ below is meant to measure the noise level at the
leaves.
\begin{definition}[Correlation at the root]\label{def:CFN_t_e}
Let $T=(V,E)$ be a tree rooted at $\rho$ with leaf set
$\partial T$.
For functions $\tilde\theta : E \to [0,1]$ and $\tilde\eta : \partial T \to [0,1]$,
let $\mathrm{CFN}(\tilde\theta,\tilde\eta)$
be the CFN model on $T$ where
$\theta(e) = \tilde\theta(e)$ for all $e$ which are not adjacent to $\partial T$, and
$\theta(e) = \tilde\theta(e) \tilde\eta(v)$ for all $e=(u,v)$,
with $v \in \partial T$.
Let
{
\[
\MajCorr(\tilde\theta,\tilde\eta) = \E\left[+\Maj(\sigma_{\partial T})\,\Big|\, \sigma_{\rho} = +1\right] =
\E\left[-\Maj(\sigma_{\partial T})\,\Big|\, \sigma_{\rho} = -1\right],
\]
}
where $\sigma$ is one sample drawn from $\mathrm{CFN}(\tilde\theta,\tilde\eta)$.
\end{definition}
\begin{proposition}[Noise reduction of majority~\cite{Mossel:04a}]\label{thm:maj_good}
Let $b$ and $\theta_{\min}$ be such that $b \theta_{\min}^2 > h^2
> 1$. Then there exist $\ell=\ell(b,\theta_{\min})$,
$\alpha=\alpha(b,\theta_{\min}) > h^{\ell} > 1$ and $\beta=\beta(b,\theta_{\min}) >
0$, such that any $\mathrm{CFN}(\theta,\eta)$ model on the $\ell$-level
$b$-ary tree satisfying $\min_{e \in E} \theta(e) \geq \theta_{\min}$ and $\min_{v \in \partial T}
\eta(v) \geq \eta_{min} >0$ also satisfies:
\begin{equation} \label{eq:maj_prop2}
\MajCorr(\theta,\eta) \geq \min\{\alpha \eta_{\min}, \beta\}.
\end{equation}
In particular if $\eta_{\min} \geq \beta$ then:
\begin{equation} \label{eq:maj_prop}
\MajCorr(\theta,\eta) \geq \beta.
\end{equation}
\end{proposition}

\paragraph{General Trees.} Recursive application of Proposition~\ref{thm:maj_good} allows the
reconstruction of the root value on any {\em balanced}
binary tree with correlation at least $\beta$. However, below we
consider {\em general} trees. In particular, when estimating the
sequence at an internal node $u$ of the phylogenetic tree, we wish to
apply Proposition~\ref{thm:maj_good} to a subtree rooted at $u$ and
this subtree need not be balanced. This can be addressed by
``completing'' the subtree into a balanced tree (with number of levels a multiple of $\ell$)
and taking
all added edges to have length $d(e) = 0$, that is, $\theta(e) = 1$.
Fix $\ell$ and $\beta$ so as to satisfy Proposition~\ref{thm:maj_good}
with $b=2$ and $\theta_{\min} = \theta$. (Recall that for the proof of
Theorem~\ref{thm} we assume $\theta(e) \geq \theta >
\theta^{\ast}$, for all edges $e$.) Consider the following
recursive function of $x = (x_1,x_2,\ldots)$:
\begin{eqnarray*}
&&\Maj_{\ell}^0(x_1) = x_1,\\
&&\Maj_{\ell}^j(x_1,\ldots,x_{2^{j{\ell}}}) = \Maj\left(
\Maj_{\ell}^{j-1}(x_1,\ldots,x_{2^{(j-1){\ell}}}),
\ldots,
\Maj_{\ell}^{j-1}(x_{2^{j{\ell}} - 2^{(j-1){\ell}} + 1},\ldots,x_{2^{j{\ell}}}),
\right),
\end{eqnarray*}
for all $j \geq 1$. Now, if $\sigma_{\partial T}$ is a character at the leaves of $T$, let us define the function $\anc{\rho}{T}(\sigma_{\partial T})$ that estimates the ancestral state at the root $\rho$ of $T$ using recursive majority as follows: 
\begin{enumerate}
\item Let $\widetilde T$ be the tree $T$ (minimally) completed with edges of $\theta$-value 1 so that
$\widetilde T$ is a complete binary tree with a number of levels a multiple of $\ell$, say $J \ell$;
\item Assign to the leaves of $\widetilde T$ the value of their ancestor leaf in $T$ under $\sigma_{\partial T}$;
\item Let $\tilde\sigma$ be the leaf states of $\widetilde T$ arranged in pre-order;
\item Compute
\begin{equation*}
\anc{\rho}{T}(\sigma_{\partial T}) := \Maj_{\ell}^J(\tilde\sigma).
\end{equation*}
\end{enumerate}
From Proposition~\ref{thm:maj_good}, we get:
\begin{proposition}[Recursive Majority]\label{cor:maj}
Let $T=(V,E)$ be a tree rooted at $\rho$ with leaf set
$\partial T$. Let $\sigma_{\partial T}$ be one sample drawn from
$\mathrm{CFN}(\theta,\eta)$ with $\theta(e) \geq \theta >
\theta^{\ast}$ for all edges $e$ and $\eta(v) = 1$ for all $v \in \partial T$.
Then, we can choose $\ell$ and $\beta > 0$ so that
\begin{equation} \label{eq:seqrecbd}
\prob\left[\anc{\rho}{T}(\sigma_{\partial T}) = \sigma_{\rho} \,\bigg|\, \sigma_{\rho} = +1\right]
= \prob\left[\anc{\rho}{T}(\sigma_{\partial T}) = \sigma_{\rho} \,\bigg|\, \sigma_{\rho} = -1\right] \geq
\frac{1+\beta}{2}.
\end{equation}
\end{proposition}
In the remainder of the paper, we will use Proposition~\ref{cor:maj} with $\theta(e) = e^{-2 d(e)}$. Also, if $\sigma_{\partial T} = \left(\sigma^i_{\partial T}\right)_{i=1}^k$ is a collection of $k$ characters at the leaves of $T$, we extend the function $\anc{\rho}{T}$ to collections of characters in the natural way:
\begin{equation} \label{eq:seqrec}
\anc{\rho}{T}(\sigma_{\partial T}) := \left(\anc{\rho}{T}(\sigma^i_{\partial T})\right)_{i=1}^k.
\end{equation}

%{
%\begin{equation} \label{eq:seqrecbd}
%\P[\textsc{SeqRec}(T,\rho)(\sigma_{\partial T}) = \sigma_{\rho} | \sigma_{\rho}] \geq
%\frac{1+\beta}{2}.
%\end{equation}
%}

%The previous theorem allows to reconstruct the root s
%balanced binary tree given values at the leaves. The estimate
%is guaranteed to have a positive correlation with the true value.
%For a balanced tree that contains
%more than $l$ levels,
%the theorem can be applied recursively and, this way, one can estimate
%internal states deep inside the tree with correlation at least
%$\beta$.

\section{Distance Estimation}\label{sec:stattech}

Throughout this section, we fix a tree $T$ on $n$ leaves.
Recall the definition of our path metric $d$ from Section~\ref{sec:definitions}.
We assume that $\sigma_{\partial T} = \left(\sigma^t_{\partial T}\right)_{t=1}^k$ are
$k$ i.i.d. samples (or characters) at the leaves
of $T$ generated by the CFN model with parameters $d(e) \leq g$ for all
$e \in \ecal(T)$ (a lower bound on $d$ is not required in this section,
but it will be in the next one).
We think of $\sigma_{\partial T} = \left((\sigma_{u}^t)_{u\in \partial T}\right)_{t=1}^k$
as $n$ sequences of length $k$ and we sometimes refer to $\sigma_u = (\sigma_u^t)_{t=1}^k$
as the ``sequence at $u$''.

\paragraph{Distance between leaves.}  As explained in Section~\ref{sec:definitions},
a basic concept in phylogenetic reconstruction is the notion of a metric on the leaves of the tree.
The distance between two leaves is a measure of the correlation between
their respective sequences. We begin with the definition of a natural distance estimator.
\begin{definition}[Correlation between sequences]
Let $u,v$ be leaves of $T$.
The following quantity
\begin{equation}\label{eq:dist}
\distance{\sigma_u}{\sigma_v} =
-\frac{1}{2}\ln\left[\left(\frac{1}{k} \sum_{t=1}^k \sigma_u^t \sigma_v^t \right)_{+}\right],
\end{equation}
is an estimate of $d(u,v)$. In cases where the sum inside the $\log$
is non-positive, we define $\distance{\sigma_u}{\sigma_v} = +\infty$.
\end{definition}

The next proposition provides a guarantee on the performance of
$\logcorr$.
The proof follows from standard concentration inequalities.
See, e.g.,~\cite{ErStSzWa:99a}.
An important point to note is that, if we use {\em short} sequences of length $O(\log n)$, the guarantee applies only to short distances, that is, distances of order $O(1)$.
\begin{proposition}[Accuracy of $\logcorr$]\label{lem:statisticsPair}
For all $\eps > 0$, $M>0$, there exists $c =
c(\eps,M; \gamma) > 0$, such that if the following conditions hold:
\begin{itemize}
\item $\mathrm{[Short\ Distance]}$ $d(u,v)<M$,
\item $\mathrm{[Sequence\ Length]}$ $k = c' \log{n}$, for some $c' > c$,
\end{itemize}
then
\begin{equation*}
\left|d(u,v)-\logcorr(\sigma_u,\sigma_v)\right|<\eps,
\end{equation*}
with probability at least $1-O(n^{-\gamma})$.
\end{proposition}
\begin{proof}
By Azuma's inequality we get
\begin{eqnarray*}
\prob\left[\logcorr(\sigma_u,\sigma_v) > d(u,v) + \eps\right]
&=& \prob\left[\frac{1}{k} \sum_{t=1}^k \sigma_u^t \sigma_v^t < e^{-2 d(u,v) - 2\eps}\right]\\
&=& \prob\left[\frac{1}{k} \sum_{t=1}^k \sigma_u^t \sigma_v^t < e^{-2 d(u,v)} - (1 - e^{-2\eps})e^{-2 d(u,v)}\right]\\
&\leq& \prob\left[\frac{1}{k} \sum_{t=1}^k \sigma_u^t \sigma_v^t < \expec[\sigma_u^1 \sigma_v^1] - (1 - e^{-2\eps})e^{-2 M}\right]\\
&\leq& \exp\left(-\frac{\left((1 - e^{-2\eps})e^{-2 M}\right)^2}{2k(2/k)^2}\right)\\
&=& n^{-c' K},
\end{eqnarray*}
with $K$ depending on $M$ and $\eps$. Above, we used that
\begin{equation*}
d(u,v) = -\frac{1}{2}\ln \expec[\sigma_u^1 \sigma_v^1].
\end{equation*}
A similar inequality holds for the other direction.
\end{proof}

\paragraph{Distance Between Internal Nodes.}
Our reconstruction algorithm actually requires us
to compute distances between {\em internal} nodes.
Note that, in general, we do not know the true sequences at
internal nodes of the tree. Therefore, we need to estimate
distances by applying (\ref{eq:dist}) to {\em reconstructed} sequences.
An obvious issue
with this idea is that reconstructed sequences are subject to
a systematic bias.
\begin{definition}[Bias]\label{def:bias}
Suppose $v$ is the root of a subtree $T_v$ of $T$. Let
%\begin{equation*}
$\widehat{\sigma}_{v} = \anc{v}{T_v}(\sigma_{\partial T_v})$.
%\end{equation*}
Then, the
quantity
\[
\lambda(T_v,v) = -\frac{1}{2}\ln(\expec[\sigma^1_v \cdot \widehat{\sigma}^1_v]_+),
\]
is called
the \emph{reconstruction bias} at $v$ on $T_v$. (Note the similarity to
(\ref{eq:dist}).) We denote
\[
\bias(\tilde g) = \sup \lambda(T_v,v),
\]
where the supremum is taken over all rooted binary trees with edge lengths
at most $\tilde g$. Note that for all $\tilde g$, $0<\tilde g < g^*$, we have by Proposition~\ref{cor:maj}
\begin{equation*}
0 < B(\tilde g) < +\infty.
\end{equation*}
We also denote
\begin{equation*}
\beta(\tilde g) = e^{-2 B(\tilde g)}.
\end{equation*}
\end{definition}
%{\bf Elchanan: Changed here}

\begin{figure}
\begin{center}
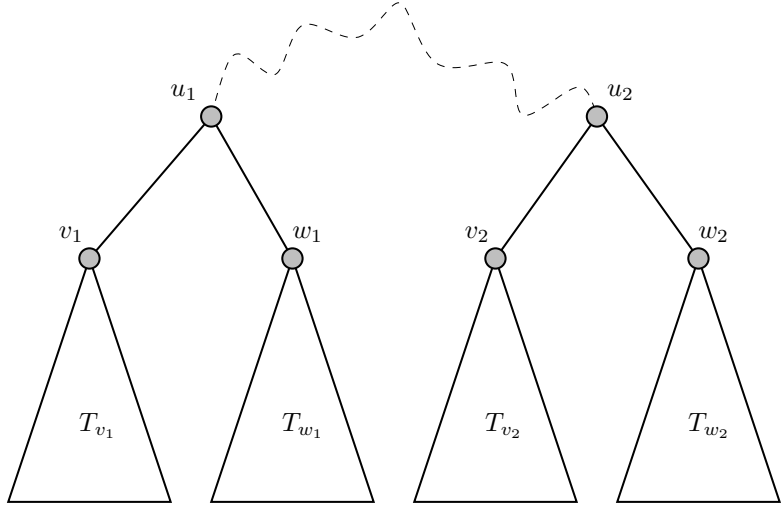\caption{Basic disjoint setup for estimator $\hatd$ and
routines \issplit\ in Figure~\ref{routine:issplit} and\ \isshort\ in Figure~\ref{routine:isshort}.
The subtrees $T_{v_1}, T_{w_1}, T_{v_2}, T_{w_2}$ are edge-disjoint
(see Definition~\ref{def:disjoint} below).}\label{fig:issplit}
\end{center}
\end{figure}

\begin{figure}
\begin{center}
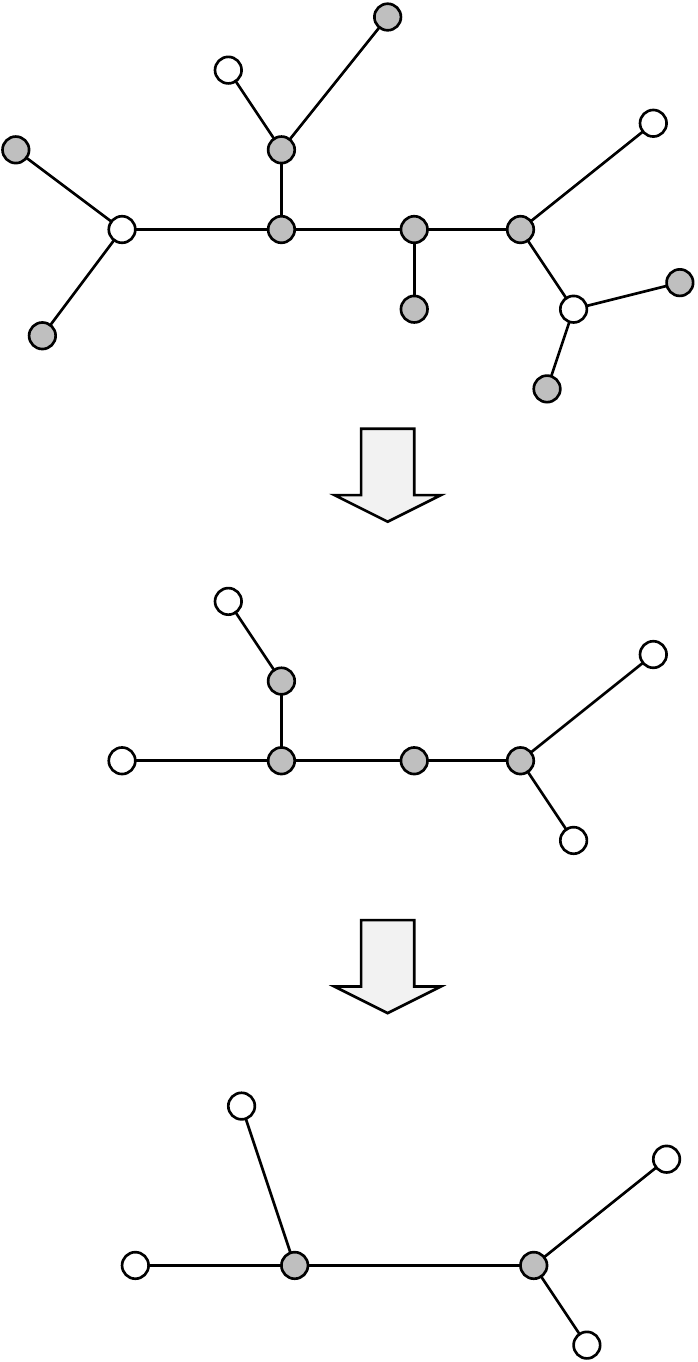\caption{Restricting the top tree to its white nodes.}\label{fig:restricted}
\end{center}
\end{figure}

\begin{figure}
\begin{center}
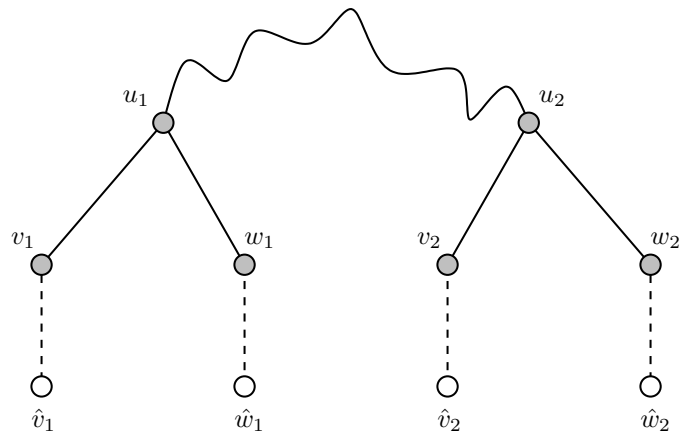\caption{Basic setup for estimator $\hatd$ where
reconstruction bias is represented by dashed edges.}\label{fig:bias}
\end{center}
\end{figure}

Since $B(g) > 0$, we cannot estimate exactly the internal sequences. Hence
Proposition~\ref{lem:statisticsPair} cannot be used directly to estimate
distances between internal nodes.
We deal with this issue as follows.
Consider the configuration in Figure~\ref{fig:issplit}. More precisely:
\begin{definition}[Basic Disjoint Setup (Preliminary Version)]
Root $T$ at an arbitrary vertex. Note that, by reversibility,
the CFN model on a tree $T$ can be rerooted arbitrarily without
changing the distribution at the leaves.
Denote by $T_x$ the subtree of $T$ rooted at $x$.
We consider two internal nodes
$u_1$, $u_2$ that are not descendants of each other in $T$.
Denote by $v_\onetwo, w_\onetwo$ the children of $u_\onetwo$
in $T_{u_\onetwo}$,
$\onetwo=1,2$. We call this configuration the \emph{Basic Disjoint Setup
(Preliminary Version)}.
\end{definition}
We need a few basic combinatorial definitions.
\begin{definition}[Restricted Subtree]
Let $V' \subseteq V$ be a subset of the vertices of $T$.
The {\em subtree of $T$ restricted to $V'$} is the tree $T'$
obtained by 1) keeping only nodes and edges on paths between
vertices in $V'$ and 2) by then contracting all paths composed
of vertices of degree 2, except the nodes in $V'$. We sometimes use the notation
$T|_{V'}:=T' $. See Figure~\ref{fig:restricted} for an example.
%(See also \cite{Mossel:07}.)
\end{definition}
\begin{definition}[Quartet]
A {\em quartet} is a set of four leaves. More generally, for any set
$\quartet = \{v_1,w_1,v_2,w_2\}$ of four nodes of $T$, we think of
$\quartet$ as a quartet on $T|_\quartet$. We say that a quartet
$\quartet = \{v_1,w_1,v_2,w_2\}$ is nondegenerate
if none of the nodes in $\quartet$ are on a path between two other nodes in $\quartet$.
There are three possible (leaf-labeled) topologies
on a nondegenerate quartet,
called {\em quartet splits}, one for each partition of the four leaves into two
pairs.
In Figure~\ref{fig:issplit}, the correct quartet split on $\quartet$ is $\{\{v_1,w_1\}\{v_2,w_2\}\}$
which we denote by $v_1 w_1|v_2 w_2$.
\end{definition}
To compute the distance between $u_1$ and
$u_2$, we think of the path between $u_1$ and $u_2$ as the
internal edge of the quartet $\quartet = \{v_1,w_1,v_2,w_2\}$.
The reconstructed sequences at $x \in \quartet$
also suffer from a systematic bias. However, we prove in Proposition~\ref{lem:statisticsQuartetPre} that the bias does not
affect the computation of the length of the \emph{internal} edge of the
quartet.
Indeed, as depicted in Figure~\ref{fig:bias}, we treat the systematic errors introduced
by the reconstruction
as ``extra edges'' attached to the
nodes in $\quartet  = \{v_1,w_1,v_2,w_2\}$, with corresponding
``endpoints'' $\hatquartet = \{\hat{v}_1, \hat{w}_1, \hat{v}_2,
 \hat{w}_2\}$.
Our estimator $\hatd$ is then obtained from a classical
distance-based computation---known in phylogenetics as the {\em
Four-Point Method}~\cite{Buneman:71}---applied to the ``extra nodes''
$\hatquartet$ (see (\ref{eq:defd}) below).
For this idea to work, it is crucial that the systematic error in the reconstructed sequences
at
$x \in \quartet$
be {\em independent} given the
true sequences at $x \in \quartet$.
This is the case when $u_1$ and $u_2$ are not
descendant of each other because of the Markov property.

We now define our distance estimator
$\hatd$. We let $\widehat{\sigma}_{x} = \anc{x}{T_x}(\sigma_{\partial T_x})$
be the reconstructed sequence at $x$.
\begin{definition}[Distance Estimator $\hatd$] \label{def:d}\text{}
Consider the Basic Disjoint Setup (Preliminary Version).
Then, we let
{\small
\begin{equation} \label{eq:defd}
\internal{\widehat{\sigma}_{v_1}}{\widehat{\sigma}_{w_1}}{\widehat{\sigma}_{v_2}}{\widehat{\sigma}_{w_2}} = \frac{1}{2} \left(
\logcorr(\widehat{\sigma}_{v_1},\widehat{\sigma}_{v_2}) +
\logcorr(\widehat{\sigma}_{w_1},\widehat{\sigma}_{w_2}) -
\logcorr(\widehat{\sigma}_{v_1},\widehat{\sigma}_{w_1}) -
\logcorr(\widehat{\sigma}_{v_2},\widehat{\sigma}_{w_2})\right).
\end{equation}
}
If one of the $\logcorr$ quantities is
$+\infty$, the function $\hatd$ is set to $+\infty$.
\end{definition}

The next proposition provides a guarantee on the performance of $\hatd$.
As in Proposition~\ref{lem:statisticsPair}, the guarantee applies only
to $O(1)$ distances, if we use sequences of length $O(\log n)$.
\begin{proposition}[Accuracy of $\hatd$]\label{lem:statisticsQuartetPre}
Consider the Basic Disjoint Setup (Preliminary Version).
Let $\quartet = \{v_1,w_1, v_2, w_2\}$. For all
$\eps > 0$, $M > 0$, there exists
$c = c(\eps, M; \gamma, g) > 0$ such that, if the following
hold:
\begin{itemize}
\item \itemname{Edge\ Lengths} For all $e \in \ecal(T)$, $d(e) \leq g < g^*$;

\item $\mathrm{[Short\ Distances]}$ We have $d(x,y)<M$, $\forall x,y \in \quartet$;

\item $\mathrm{[Sequence\ Length]}$ The sequences used have length $k = c' \log{n}$, for some $c'
> c$;
\end{itemize}
then, with probability at least $1-O(n^{-\gamma})$,
\begin{equation*}
\left|d(u_1,u_2)-\internal{\widehat{\sigma}_{v_1}}{\widehat{\sigma}_{w_1}}{\widehat{\sigma}_{v_2}}{\widehat{\sigma}_{w_2}}\right|<\eps.
\end{equation*}
\end{proposition}
\begin{proof}
First note that, for all $t\in\{1,\ldots,k\}$, conditioned on
$\{\sigma_{x}^t\}_{x\in\quartet}$,
the values $\{\recseq{x}^t\}_{x\in\quartet}$ are independent by the Markov property.
Therefore, as we observed before, the systematic errors introduced by the reconstruction process can be treated
as ``extra edges'', with endpoints
$\hatquartet = \{\hat{v}_1, \hat{w}_1, \hat{v}_2, \hat{w}_2\}$ (see Figure~\ref{fig:bias}).
Moreover, from (\ref{eq:seqrecbd}), it follows that, $\forall x \in
\quartet$, $\forall t\in\{1,\ldots,k\}$,
\begin{equation}\label{eq:seqrecbd2}
\prob\left[\widehat{\sigma}_x^t = \sigma_x^t| \sigma_x^t\right]  \ge \frac{1+\beta(g)}{2},
\end{equation}
where $\beta(g)>0$ is specified by Definition~\ref{def:bias}.
From (\ref{eq:seqrecbd2}),
the lengths of the ``extra edges''---equal to the reconstruction
biases---are at most $\bias(g)$.

Appealing to Proposition~\ref{lem:statisticsPair} with
constants $\frac{\eps}{2}$ and $M + 2\bias(g)$
and the union bound, we argue that with probability
at least $1-O(n^{-\gamma})$ for all $x,y \in \quartet$
\begin{equation*}
\left|\logcorr(\widehat{\sigma}_{x},\widehat{\sigma}_{y}) - d(\hat{x},\hat{y})\right| < \frac{\eps}{2},
\end{equation*}
where we extend $d$ to the ``extra nodes'' $\hat{x}, \hat{y}$ by setting
\begin{equation*}
d(x, \hat{x}) = -\frac{1}{2}\ln(\expec[\sigma_x \widehat{\sigma}_x]_+),
\end{equation*}
and likewise for $\hat{y}$ (see Figure~\ref{fig:bias}).
Observing that
$d(u_1,u_2) = \frac{1}{2}\left( d(\hat{v}_1,\hat{v}_2)+d(\hat{w}_1,\hat{w}_2) - d(\hat{v}_1,\hat{w}_1)- d(\hat{v}_2,\hat{w}_2)\right)$ and using the above we get
\begin{equation*}
\left|d(u_1,u_2)-\internal{\widehat{\sigma}_{v_1}}{\widehat{\sigma}_{w_1}}{\widehat{\sigma}_{v_2}}{\widehat{\sigma}_{w_2}}\right|<\eps.
\end{equation*}
%}
\end{proof}

\paragraph{Distances Between Restricted Subtrees.}
In fact, we need to apply Proposition~\ref{lem:statisticsQuartetPre}
to {\em restricted subtrees} of the true tree.
Indeed, in the reconstruction algorithm, we maintain a ``partially reconstructed
subforest'' of the true tree, that is, a collection
of restricted subtrees of $T$.
%\begin{definition}[Restricted Subforest]
%The forest $\forestno = \{T_1, \ldots, T_\alpha\}$
%is a {\em restricted subforest of $T$} if there is a collection
%of disjoint subsets of leaves $L_1, \ldots, L_\alpha$ of $T$
%,
%$L_\onetwo \subseteq \partial T$, $\onetwo = 1,\ldots, \alpha$,
%with $L_\alpha \cap L_\beta = \emptyset$ for $\alpha \neq \beta$,
%such that $\forestno = \{T|_{L_1},\ldots,T|_{L_\alpha}\}$.
%We say that $\forestno$ is a {\em rooted restricted subforest}
%if each $T_\onetwo$, $\onetwo = 1,\ldots, \alpha$, is rooted.
%We denote by $\rho(\forestno)$ the set of roots of $\forestno$.
%(See also \cite{Mossel:07}.)
%\end{definition}
For this more general setup,
we use the routine \distest\ detailed in Figure~\ref{fig:distest}.
%The routine takes as input a rooted restricted subforest $\fcal$ of the true tree $T$
%and the nodes in $\forestno$.
%The forest $\fcal$ is typically a collection of rooted induced subtrees of $T$
%with edges pointing away from the roots. We denote
%by $T_x^\fcal$ the subtree of $\fcal$ rooted at $x$.
\begin{figure*}[h]
\framebox{
\begin{minipage}{16cm}
{\small \textbf{Algorithm} \distest\\
\textit{Input:} Two nodes $u_1, u_2$;
a rooted forest $\fcal$;
accuracy radius $\erada$;\\
\textit{Output:} Distance estimate $\nu$;
\begin{itemize}
\item \itemname{Children} Let
$w_\onetwo, v_\onetwo$ be the children of $u_\onetwo$ in $\fcal$
for $\onetwo = 1,2$ (if $u_\onetwo$ is a leaf,
set $w_\onetwo = v_\onetwo = u_\onetwo$);

\item $\mathrm{[Sequence\ Reconstruction]}$ For $x \in \{v_1, w_1, v_2, w_2\}$, set
$\recseq{x} = \anc{x}{T_x^{\fcal}}(\sigma_{\partial T_x^{\fcal}})$;

\item \itemname{Accuracy\ Cutoff}
If there is $\{x,y\} \subseteq \{v_1, w_1, v_2, w_2\}$ such that
\begin{equation*}
\distance{\recseq{x}}{\recseq{y}} > \erada,
\end{equation*}
return $+\infty$;

\item $\mathrm{[Distance\ Computation]}$ Otherwise return
\begin{equation*}
\nu = \internal{\recseq{v_1}}{\recseq{w_1}}{\recseq{v_2}}{\recseq{w_2}}.
\end{equation*}
\end{itemize}
}
\end{minipage}
} \caption{Routine \distest.} \label{fig:distest}
\end{figure*}

\begin{figure}
\begin{center}
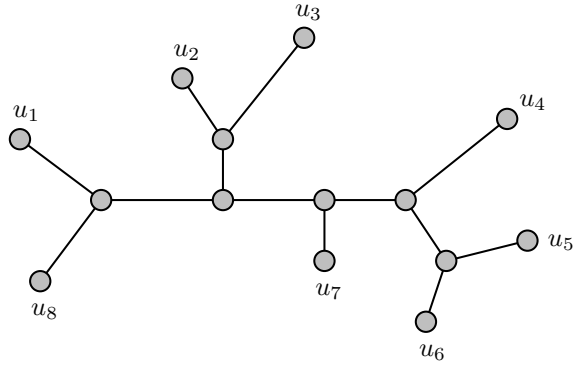\caption{
The subtrees $T|_{\{u_1,u_2,u_3,u_8\}}$
and $T|_{\{u_4,u_5,u_6,u_7\}}$ are edge-disjoint.
The subtrees $T|_{\{u_1,u_5,u_6,u_8\}}$
and $T|_{\{u_2,u_3,u_4,u_7\}}$ are edge-sharing.
}\label{fig:sharing}
\end{center}
\end{figure}

\begin{figure}
\begin{center}
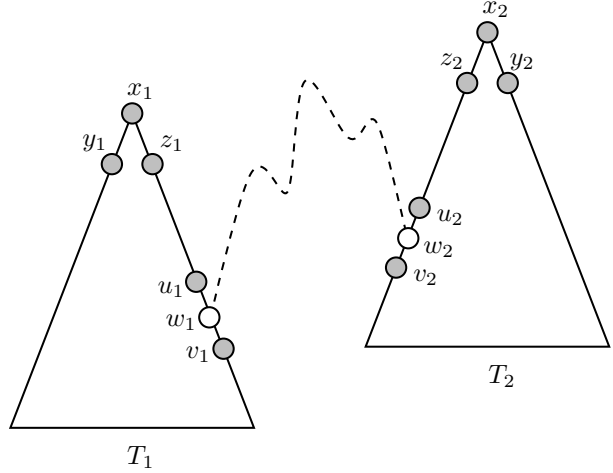\caption{
Basic Disjoint Setup (General).
The rooted subtrees $T_{1}, T_{2}$ are edge-disjoint
but are not assumed to be dangling.
The white nodes may not be in the {\em restricted} subtrees
$T_{1}, T_{2}$.
The case $w_1 = x_1$ and/or $w_2 = x_2$ is possible.
Note that if we root the tree at any node along the dashed path, the
subtrees rooted at $y_1$ and $y_2$ are edge-disjoint and dangling
(unlike $T_1$ and $T_2$).}\label{fig:remote}
\end{center}
\end{figure}

To generalize Proposition~\ref{lem:statisticsQuartetPre}, we need a few
definitions.
First, the notion of edge disjointness is borrowed from~\cite{Mossel:07}.
\begin{definition}[Edge Disjointness] \label{def:disjoint}
Denote by $\path_T(x,y)$ the path (sequence of edges) connecting
$x$ to $y$ in $T$. We say that two restricted subtrees $T_1, T_2$ of $T$ are
{\em edge disjoint} if
\[
\path_T(x_1,y_1) \cap \path_T(x_2,y_2) = \emptyset,
\]
for all $x_1, y_1 \in \lcal(T_1)$ and $x_2, y_2 \in \lcal(T_2)$. We
say that $T_1,T_2$ are {\em edge sharing} if they are not edge
disjoint. See Figure~\ref{fig:sharing} for an example. (If $T_1$ and $T_2$ are directed, we take this
definition to refer to their underlying undirected version.)
\end{definition}
%{\bf\red [SR: Expanded definition below.]}
\begin{definition}[Legal Subforest]
We say that a tree is {\em a rooted full binary tree} if all its internal
nodes have degree 3 except for the root which has degree 2. A restricted
subtree $T_1$ of $T$ is a {\em legal subtree} of $T$ if:
\begin{enumerate}
\item It is
a rooted full binary tree, that is, if it has a unique node of degree 2---its root;
\item And
if all its leaves are leaves of $T$.
\end{enumerate}
Moreover, we say that a forest
\begin{equation*}
\forestno = \{T_1,\ldots,T_\alpha\},
\end{equation*}
is {\em legal subforest} of $T$ if the $T_\onetwo$'s are
\emph{edge-disjoint} legal subtrees of $T$.
We denote by $\rho(\forestno)$ the set of roots of $\forestno$.
\end{definition}
\begin{definition}[Dangling Subtrees]
%Let $\forestno$ be a rooted restricted subforest of $T$.
We say that two edge-disjoint legal subtrees $T_1$, $T_2$ of $T$
are {\em dangling} if there is a choice of root $\rho^*$ for $T$
\emph{not in $T_1$ or $T_2$}
that is consistent with the rooting of both $T_1$ and $T_2$,
that is, pointing the edges of $T$ away from $\rho^*$ is such
that the edges of $T_1$ and $T_2$ are directed away from their respective roots.
\end{definition}
See Figure~\ref{fig:remote} below for an example where two legal, edge-disjoint subtrees
are \emph{not} dangling.
We generalize the basic configuration of Figure~\ref{fig:issplit} as follows:
\begin{definition}[Basic Disjoint Setup (Dangling)]
Let $T_1 = T_{u_1}$ and $T_2 = T_{u_2}$ be two legal subtrees of $T$ rooted
at $u_1$ and $u_2$ respectively.
Assume further that
$T_1$ and $T_2$ are {\em edge-disjoint} and {\em dangling}.
Denote by $v_\onetwo, w_\onetwo$ the children of $u_\onetwo$
in $T_\onetwo$,
$\onetwo=1,2$. If $u_\onetwo$ is a leaf, we let instead $v_\onetwo = w_\onetwo = u_\onetwo$.
We call this configuration the \emph{Basic Disjoint Setup
(Dangling)}.
\end{definition}

We make another important remark about the more general setup considered here.
Note that even though
the true tree satisfies the ``short edge length'' condition (that is, $d(e)\leq g$)
by assumption, partially reconstructed subtrees {\em may not}---because
their edges are actually {\em paths} in the true tree.
Therefore, an important step of the reconstruction algorithm
is to make sure that all restricted edges
of the partially reconstructed subforest are short enough for recursive majority
to be accurate.
(This explains why we need to use two constants $g < g'$ below $g^*$.)
Refer to routine \cherryid\ in Figure~\ref{routine:cherryid} in Section~\ref{sec:algvars}.

We then obtain the following generalization of
Proposition~\ref{lem:statisticsQuartetPre}.
Note that the routine \distest\ in Figure~\ref{fig:distest}
uses an ``accuracy cutoff'' $\erada$.
As we show in the proof below, this ensures that distances that are too long
are rejected, in which case the value $+\infty$ is returned.
\begin{proposition}[Accuracy of \distest]\label{lem:statisticsQuartet}
Consider the Basic Disjoint Setup (Dangling).
Let $\quartet = \{v_1,w_1, v_2, w_2\}$. For all
$\eps > 0$, $M > 0$, and $\erada > M + 2\bias(g')$, there exists
$c = c(\eps, M, \erada; \gamma, g') > 0$ such that, if the following
hold:
\begin{itemize}
\item $\mathrm{[Edge\ Length]}$ It holds that $d(e)\leq g' < g^*$, $\forall e \in \ecal(T_{x})$, $x\in\quartet$;

\item $\mathrm{[Sequence\ Length]}$ The sequences used have length $k = c' \log{n}$, for some $c'
> c$;
\end{itemize}
then, with probability at least $1-O(n^{-\gamma})$, the following holds:
letting $\nu$ be the output of \distest\ in Figure~\ref{fig:distest}, we have that
if one of the following hold
\begin{enumerate}
\item $\mathrm{[Short\ Distances]}$ We have $d(x,y)<M$, $\forall x,y \in \quartet$;
\item $\mathrm{[Finite\ Estimate]}$ $\nu < +\infty$;
\end{enumerate}
then 
\begin{equation*}
\left|d(u_1,u_2)-\nu\right|<\eps.
\end{equation*}
\end{proposition}
\begin{proof}
The first half of the proposition follows immediately from Propositions~\ref{lem:statisticsPair} and~\ref{lem:statisticsQuartetPre}. Refer to Figure \ref{fig:bias}. 

The second part follows from a ``double window'' argument as in~\cite[Theorem 9]{ErStSzWa:99a}.
Let $0 < \tau < {1\over 2}$ such that
\begin{equation*}
\erada = M + 2B(g') - \frac{1}{2}\ln (2 \tau).
\end{equation*}
By assumption, $\nu < +\infty$ and therefore for all pairs
$\{x,y\} \subseteq \{v_1, w_1, v_2, w_2\}$ we have
\begin{equation}\label{eq:fourtau}
\distance{\recseq{x}}{\recseq{y}} \leq M + 2B(g')  - \frac{1}{2}\ln (2 \tau).
\end{equation}
It follows directly from Azuma's inequality that a pair $\{x,y\}$ such that
\begin{equation*}
d(\hat{x},\hat{y}) \ge M + 2B(g') - \frac{1}{2}\ln \tau,
\end{equation*}
satisfies (\ref{eq:fourtau}) with probability at most
\begin{eqnarray*}
\prob\left[\distance{\recseq{x}}{\recseq{y}} \leq M + 2B(g') - \frac{1}{2}\ln 2 \tau\right]
&=& \prob\left[\frac{1}{k} \sum_{t=1}^k \recseq{x}^t \recseq{y}^t \geq 2\tau e^{-2[M + 2B(g')]}\right]\\
&\leq& \prob\left[\frac{1}{k} \sum_{t=1}^k \recseq{x}^t \recseq{y}^t \geq e^{-2 d(\hat{x},\hat{y})} + \tau e^{-2[M + 2B(g')]}\right]\\
&\leq& \exp\left(-\frac{\left(\tau e^{-2[M + 2B(g')]}\right)^2}{2k(2/k)^2}\right)\\
&=& n^{-c' K},
\end{eqnarray*}
for a constant $K > 0$ depending on $M, g', \tau$.
The result follows by applying the first part of the proposition with $M$
replaced by $M + 2B(g') - \frac{1}{2}\ln \tau$.
(Note that $d(x,y) \le d(\hat{x},\hat{y})$, for all $x,y \in \quartet$.)
\end{proof}

We also need a simpler variant of \distest\ (Figure~\ref{fig:distest}) whose purpose is to test
whether the internal path of a quartet is longer than $g$.
We record this variant in Figure~\ref{routine:isshort}.
\begin{figure*}[h]
\framebox{
\begin{minipage}{16cm}
{\small \textbf{Algorithm} \isshort\\
\textit{Input:} Two pairs of nodes $(v_1, w_1)$, $(v_2, w_2)$;
rooted forest $\forestno$;
accuracy radius $\erada$;
tolerance $\tol$;\\
\textit{Output:} Boolean value and length estimate;
\begin{itemize}
\item $\mathrm{[Sequence\ Reconstruction]}$ For $x \in \{v_1, w_1, v_2, w_2\}$, set
$\recseq{x} = \anc{x}{T_x^{\fcal}}(\sigma_{\partial T_x^{\fcal}})$;

\item \itemname{Accuracy\ Cutoff}
If there is $\{x,y\} \subseteq \{v_1, w_1, v_2, w_2\}$ such that
\begin{equation*}
\distance{\recseq{x}}{\recseq{y}} > \erada,
\end{equation*}
return $+\infty$;

\item \itemname{Internal\ Edge\ Length}
Set
\begin{equation*}
\nu = \internal{\recseq{v_1}}{\recseq{w_1}}{\recseq{v_2}}{\recseq{w_2}};
\end{equation*}
\item \itemname{Test} If $\nu < g + \tol$ return
$(\mathrm{\true}, \nu)$, o.w. return $(\mathrm{\false},0)$;
\end{itemize}
}
\end{minipage}
} \caption{Routine \isshort.} \label{routine:isshort}
\end{figure*}

\paragraph{Detecting Long Distances When $T_1$ And $T_2$ Are Not Dangling.}
%\paragraph{Remote Subtrees.}
%Finally, we need one last important property of the estimator
%$\hatd$, namely that it can detect whether two edge-disjoint subtrees are
%sufficiently far apart.
Roughly speaking, our reconstruction algorithm
works by progressively merging
subtrees that are close in the true tree. (See Section~\ref{sec:algvars} for
further details.) Hence, the algorithm needs to tell whether or not two subtrees are
sufficiently close to be considered for this merging operation. However,
as we explain in the next sections, we cannot guarantee that the
Basic Disjoint Setup (Dangling) in Proposition~\ref{lem:statisticsQuartet}
applies to all situations encountered during the execution of the algorithm.
Instead we use a special ``multiple test'' to detect long distances.
This test is performed by the routine \distmet\ detailed in Figure~\ref{fig:distmet}.

The routine has a further important property. During the course of the algorithm,
since we only have partial information about the structure of the tree,
we may not always know whether or not two subtrees are dangling---and therefore
whether or not \distest\ in Figure~\ref{fig:distest} returns accurate distance estimates. 
The ``multiple test''
in \distmet\ is such that, if the routine returns a {\em finite} estimate, that
estimate is accurate. We proceed with an explanation of these properties.

We further generalize the basic configuration of Figure~\ref{fig:issplit} as follows:
\begin{definition}[Basic Disjoint Setup (General)]\label{def:bdsgeneral}
Let $T_1 = T_{x_1}$ and $T_2 = T_{x_2}$ be two restricted subtrees of $T$ rooted
at $x_1$ and $x_2$ respectively.
Assume further that
$T_1$ and $T_2$ are {\em edge-disjoint}, but not necessarily {\em dangling}.
Denote by $y_\onetwo, z_\onetwo$ the children of $x_\onetwo$
in $T_\onetwo$,
$\onetwo=1,2$.
Let $w_\onetwo$ be the node in $T$ where the path
between $T_1$ and $T_2$ meets $T_\onetwo$, $\onetwo = 1,2$.
Note that $w_\onetwo$ may not be in $T_\onetwo$ since $T_\onetwo$ is {\em restricted}, $\onetwo = 1,2$.
If $w_\onetwo \neq x_\onetwo$, assume without loss of generality that $w_\onetwo$
is in the subtree of $T$ rooted at $z_\onetwo$, or on the edge $(x_\onetwo, z_\onetwo)$, $\onetwo = 1,2$.
We call this configuration the \emph{Basic Disjoint Setup
(General)}.
See
Figure~\ref{fig:remote}.
Let $d(T_1,T_2)$ be the length of the path between $w_1$ and
$w_2$ (in the path metric $d$).
\end{definition}

The key point to note is that when computing the distance
between $y_1$ and $y_2$ rather than the distance between $x_1$ and $x_2$, then
the assumptions of Proposition~\ref{lem:statisticsQuartet}
are satisfied (as the claim holds in this case by rooting the tree at any node
along the path connecting $w_1$ and $w_2$).
%{\bf Elchanan: added explanation. May also be useful to add the root for which
%this works in the figure}
Hence, if $T_1$ and $T_2$ are far apart,
%(Case~\ref{case:remote}),
the distance between $y_1$ and $y_2$ is correctly estimated as being large. On the other hand,
%in Case~\ref{case:dangling}, the assumptions of Proposition~\ref{lem:statisticsQuartet}
if $T_1$ and $T_2$ are dangling and close, then
distances between all pairs in $\{y_1, z_1\}\times\{y_2, z_2\}$ are accurately estimated
as being short.

%We consider the following situation, which we call the ``basic remote setup
%with threshold $M$'' (in both cases below, $T_1$ and $T_2$ are assumed edge-disjoint):
%\begin{enumerate}
%\item \label{case:remote} Either $d(T_1, T_2) > M$;
%\item \label{case:dangling} Or $d(T_1, T_2) \leq M$ and $T_1$ and $T_2$ are dangling.
%\end{enumerate}
%We immediately get the following proposition.
%{\bf Elchanan: $h$ is not defined below}

\begin{proposition}[Accuracy of \distmet]\label{lem:statisticsRemote}
Consider the Basic Disjoint Setup (General) with $\forestno = \{T_1, T_2\}$
and $\quartet = \{y_1,z_1,y_2,z_2\}$.
For all
$\eps > 0$, $M > 0$, and $\erada > M+2\bias(g') + 4g'$ there exists
$c = c(\eps, M, \erada; \gamma, g') > 0$
such that, if the following
hold:
\begin{itemize}
\item $\mathrm{[Edge\ Length]}$ It holds that $d(e)\leq g' < g^*$, $\forall e \in \ecal(T_{x})$, $x\in\quartet$; also, $d(x_\onetwo, y_\onetwo), d(x_\onetwo, z_\onetwo) \leq g' $, $\onetwo=1,2$;

\item $\mathrm{[Weight\ Estimates]}$ We are given weight estimates
$$h\,:\,\ecal(\forestno) \to \real_+,$$ such that $|h(e) - d(e)|<\tol/16$, $\forall e \in \ecal(\forestno)$;

\item $\mathrm{[Sequence\ Length]}$ The sequences used have length $k = c' \log{n}$, for some $c'
> c$;
\end{itemize}
then, with probability at least $1-O(n^{-\gamma})$, the following holds:
letting $\nu$ be the output of \distmet\ in Figure~\ref{fig:distmet}, we have that
if one of the following hold
\begin{enumerate}
\item \itemname{Dangling\ Case} $T_1$ and $T_2$ are dangling
and $d(T_1, T_2) < M$;
\item \itemname{Finite\ Estimate} $\nu < +\infty$;
\end{enumerate}
then
\begin{equation*}
\left|\nu - d(x_1, x_2)\right| < \tol.
\end{equation*}
\end{proposition}
\begin{proof}
Choose $c$ as in Proposition~\ref{lem:statisticsQuartet} with
parameters $\tol/16$ and $M + 4g'$. As in Figure~\ref{fig:distmet}, let
\begin{equation*}
\metricno'(r_1, r_2) := \distest(r_1, r_2; \forestno; \erada) - h(x_1, r_1) - h(x_2, r_2),
\end{equation*}
for all pairs
$(r_1, r_2) \in \{y_1, z_1\}\times\{y_2, z_2\}$ (where $\distest$ is defined in Figure~\ref{fig:distest}).
In {\em Case 1}, the result follows
directly from Proposition~\ref{lem:statisticsQuartet} and the remarks above
Proposition~\ref{lem:statisticsRemote}. In particular, note that, with
probability $1-O(n^{-\gamma})$,
\begin{equation*}
\left|\metricno'(r_1,r_2) - d(x_1,x_2)\right| < \frac{3 \tol}{16},
\end{equation*}
for all $(r_1, r_2) \in \{y_1, z_1\}\times\{y_2, z_2\}$ and
therefore
\begin{equation*}
\max\left\{\left|\metricno'(r_1^{(1)}, r_2^{(1)}) - \metricno'(r_1^{(2)}, r_2^{(2)})\right|
: (r_1^{(\onetwo)}, r_2^{(\onetwo)}) \in \{y_1, z_1\}\times\{y_2, z_2\}, \onetwo = 1,2\right\}
< 2\left(\frac{3\tol}{16}\right) < \tol/2,
\end{equation*}
and \distmet\ returns a finite value
(namely $\nu := \metricno'(z_1, z_2)$; see Figure~\ref{fig:distmet})
that is accurate within $3\tol/16 < \tol$.

In {\em Case 2}, the condition $\nu < +\infty$ implies in particular
that all four distance estimates are equal up to $\tol/2$. From the remark
above the statement of the proposition, at least one such distance,
say $\metricno'(y_1,y_2)$ w.l.o.g., is computed within the Basic Disjoint
Setup (Dangling) and we can therefore apply Proposition~\ref{lem:statisticsQuartet} again.
In particular, we have
\begin{equation*}
\left|\metricno'(y_1,y_2) - d(x_1,x_2)\right| < \frac{3 \tol}{16},
\end{equation*}
and therefore
\begin{equation*}
\left|\nu - d(x_1,x_2)\right| < \frac{\tol}{2} + \frac{3 \tol}{16} < \tol,
\end{equation*}
(where again $\nu := \metricno'(z_1, z_2)$).
\end{proof}

\paragraph{Remark.} Note that \distmet\ outputs $+\infty$ in two very distinct situations. 
In particular, if \distmet\ outputs $+\infty$, then either
the subtrees are too far to obtain an accurate estimate {\em or} $T_1$ and $T_2$ are not dangling
(or both).
This convention will turn out to be convenient in the statement
of the full algorithm.

%\begin{proof}
%This follows immediately from:
%\begin{lemma}[Correlation of Antisymmetric Functions]
%Let $T_1$ and $T_2$ be edge disjoint subtrees of $T$ whose distance
%is at least $\alpha g$.
%For $i = 1,2$, let $\phi_i : \{-1,1\}^{\vcal(T_i)} \to \{-1,1\}$
%be an antisymmetric function.
%If $\sigma$ is a character generated by the CFN model, then
%$$
%\corr(\phi_1(\sigma_{T_1}),\phi_2(\sigma_{T_2})) \leq \exp(-2 \alpha g).
%$$
%%\end{lemma}
%\begin{proof}
%We use the so-called ``random cluster''
%representation of the model (see e.g.~\cite{Grimmett:04} for more details).
%In this representation, an edge $e$ acts as
%follows:
%\begin{itemize}
%\item
%with probability $\exp(-2d(e))$ the two endpoints of the edge are
%identical,
%\item
%with probability $1-\exp(-2d(e))$ the two endpoints are independent.
%\end{itemize}
%It is now easy to see that if $r$ is the length of the path connecting
%the two subtrees, then with probability $1-e^{-2r}$ the measures on the
%two subtrees are independent. This contributes $0$ to the correlation.
%In the other case, we get a contribution of at most $1$ in absolute
%value (by symmetry, the variance is $1$).
%Thus the correlation is bounded by $e^{-2r}$ in aboslute value as needed.
%\end{proof}
%\end{proof}
\begin{figure*}[h]
\framebox{
\begin{minipage}{16cm}
{\small \textbf{Algorithm} \distmet\\
\textit{Input:} Two nodes $x_1, x_2$;
a rooted forest $\forestno$; edge lengths $\{h(e)\}_{e\in\forestno}$;
accuracy radius $\erada$; tolerance $\tol$;\\
\textit{Output:} Distance $\nu$;
\begin{itemize}
\item \itemname{Children} Let
$y_\onetwo, z_\onetwo$ be the children of $x_\onetwo$ in $\fcal$
for $\onetwo = 1,2$ (if $x_\onetwo$ is a leaf,
set $z_\onetwo = y_\onetwo = x_\onetwo$);

\item $\mathrm{[Distance\ Computations]}$ For all pairs
$(r_1, r_2) \in \{y_1, z_1\}\times\{y_2, z_2\}$, compute
\begin{equation*}
\metricno'(r_1, r_2) := \distest(r_1, r_2; \forestno; \erada) - h(x_1, r_1) - h(x_2, r_2);
\end{equation*}

\item \itemname{Multiple\ Test} If
\begin{equation*}
\max\left\{\left|\metricno'(r_1^{(1)}, r_2^{(1)}) - \metricno'(r_1^{(2)}, r_2^{(2)})\right|
: (r_1^{(\onetwo)}, r_2^{(\onetwo)}) \in \{y_1, z_1\}\times\{y_2, z_2\}, \onetwo = 1,2\right\} < \tol/2,
\end{equation*}
return $\nu := \metricno'(z_1, z_2)$,
otherwise return $\nu := +\infty$
(return $\nu := +\infty$ if any of the distances above is $+\infty$).
\end{itemize}
}
\end{minipage}
} \caption{Routine \distmet.} \label{fig:distmet}
\end{figure*}

%\section{Combinatorial Tools}\label{sec:combtech}

\section{Quartet Tests}\label{sec:combtech}

Before giving a complete description of the reconstruction algorithm,
we introduce some important combinatorial tools. As we discussed
previously, we make use of a key concept from phylogenetics---the notion
of a quartet. In the previous section, we showed how to estimate accurately
distances between internal nodes of the tree. In this section, we explain
how to use such estimates to perform topological tests on quartets.
Those tests form the basic building blocks of the combinatorial
component of the algorithm.
As before, we fix a tree $T$ on $n$ leaves. Moreover, we assume that our path
metric $d$ on $T$ satisfies $d(e) \geq f$ for all $e \in T$. Note that
an upper bound on $d$ is not explicitly required in this section since we assume
that appropriate distance estimates are handed to us.

\paragraph{Splits.}
%Our algorithm is based on quartet techniques from Phylogenetics.
%Recall that a {\em quartet} is a set of four leaves (possibly of a restricted subtree
%of the original tree).
%The basic test we perform on quartets is to decide their topology.
The routine \issplit\ in Figure~\ref{routine:issplit}
performs a classic test to decide
the correct split of a quartet. It is based on the so-called
{\em Four-Point Method}~\cite{Buneman:71}.
%See also Figure~\ref{fig:issplit}.
\begin{figure*}[h]
\framebox{
\begin{minipage}{16cm}
{\small \textbf{Algorithm} \textsc{IsSplit}\\
\textit{Input:} Two pairs of nodes $(v_1, w_1)$, $(v_2, w_2)$; a distance matrix $\metricno$
on these four nodes;\\
\textit{Output:} \true\ or \false;
\begin{itemize}
\item \itemname{Internal\ Edge\ Length}
Set
\begin{equation*}
\nu = \frac{1}{2}\left(\metricno(w_1,w_2)+\metricno(v_1,v_2)-\metricno(w_1,v_1)-\metricno(w_2,v_2)\right);
\end{equation*}
(set $\nu = +\infty$ if any of the distances is $+\infty$)
\item \itemname{Test} If $\nu < f/2$ return
\false, o.w. return \true.
\end{itemize}
}
\end{minipage}
} \caption{Routine \issplit.} \label{routine:issplit}
\end{figure*}
%\begin{figure}\centering
%\includegraphics[width=10cm, height=7cm]{FourPointGood.eps} \caption{Computing distances from
%closest common ancestor in {\sc IsSplit}.} \label{fig:fourpointgood}
%\end{figure}
Proposition~\ref{lem:issplit} guarantees the correctness of \issplit.
Its proof is omitted. We consider once again the Basic Disjoint Setup (Dangling)
of Section~\ref{sec:stattech}.
\begin{proposition}[Correctness of \issplit]\label{lem:issplit}
Consider the Basic Disjoint Setup (Dangling). Let
\begin{equation*}
\quartet=\{v_1, w_1, v_2, w_2\},
\end{equation*}
and let $\metricno$ be the
distance matrix on the four nodes of $\quartet$ passed to \issplit\ in Figure~\ref{routine:issplit}. Assume that
$d(e) \geq f$ for all edges on the subtree restricted to $\quartet$.
If
\begin{equation*}
\left|d(x,y) - \metricno(x,y)\right| < \frac{f}{4},\qquad \forall x,y \in \quartet,
\end{equation*}
then the call
\issplit$((v_1, w_1),(v_2, w_2); \metricno)$ returns \true, whereas the calls
\issplit$((v_1, w_2),(v_2, w_1); \metricno)$ and
\issplit$((v_1, v_2),(w_1, w_2); \metricno)$ return \false.
\end{proposition}

\paragraph{Collisions.}
As we discussed before, for routines \distest\ (Figure~\ref{fig:distest}) and \issplit\ (Figure~\ref{routine:issplit})
to work, we need a configuration as in Figure~\ref{fig:issplit}
where two edge-disjoint subtrees are connected by a path
that lies ``above'' them. However, for reasons that will be explained
in later sections (see also the discussion in Section~\ref{sec:stattech}),
in the course of the reconstruction algorithm we
have to deal with configurations where two subtrees
are not dangling, as in Figure~\ref{fig:remote}.
We use the following definitions. Recall the Basic Disjoint Setup (General)
in Figure~\ref{fig:remote}.
%We define next a special case of such configuration, which
%we call a {\em collision}.
%See Figure~\ref{fig:collision} for an illustration.
%\begin{figure}
%\begin{center}
%\input{collision.pstex_t}\caption{A collision. Compare to the
%basic disjoint setup in Figure~\ref{fig:issplit}.}\label{fig:collision}
%\end{center}
%\end{figure}
%Although we have not yet specified the details of our
%algorithm, we have defined its basic characteristic as that of a
%greedy algorithm that maintains a forest of disjoint rooted
%subtrees of the model tree $T$ and, at each step, uses local
%metric information to join trees of the forest and increase its
%size. Note that a property that would be useful to maintain
%along this incremental process would be that, at each step of the
%algorithm, the part of $T$ that is missing from the current forest
%``lies above'' the roots of the trees in the forest. But, since
%only partial information is known about the tree metric at each
%step of the execution, it will be unavoidable that some joins made
%by the algorithm violate this desired property and create what we
%define as {\em collisions} in Definition \ref{def:collisions}. On
%the other hand, as the forest grows new metric information is
%discovered and collisions can be identified.
\begin{definition}[Collisions]\label{def:collisions}
Suppose that $T_1$ and $T_2$ are legal subtrees of
$T$ and suppose they are {\em not} dangling.
We say that {\em $T_1$ collides into $T_2$ at edge $e_2=(u_2,v_2)$
($u_2$ is the parent of $v_2$ in $T_2$)}, if the path
\begin{equation*}
\path_T(\rho(T_1),\rho(T_2))
\end{equation*}
has non-empty intersection with
edge $e_2$ (i.e., with $\path_T(u_2,v_2)$) but with no other edge in the subtree of $T_2$ rooted at
$v_2$. See Figure~\ref{fig:remote}.
We sometimes say that the trees $T_1$, $T_2$ {\em collide}. We say that
{\em the collision is within distance $M$} if $d(T_1,T_2) \leq M$.
(See Definition~\ref{def:bdsgeneral}.)
\end{definition}
%%\begin{definition}[Collisions]
%%Suppose that $T_1$ and $T_2$ are edge disjoint rooted subtrees of
%%$T$. We say that {\em $T_1$ and $T_2$ collide at distance $d$}, if
%%the path $\path_T(\rho(T_1),\rho(T_2))$ has non-empty intersection
%%with $\ecal(T_1)\cup\ecal(T_2)$ and the length of the shortest
%%path between $T_1$ and $T_2$ is at most $d$. In other words, $T_1$
%%and $T_2$ collide at distance $d$, if the shortest path between
%%$T_1$ and $T_2$ is of length at most $d$ and this path does not
%%contain either $\rho(T_1)$ or $\rho(T_2)$.
%%\end{definition}

\begin{figure}
\begin{center}
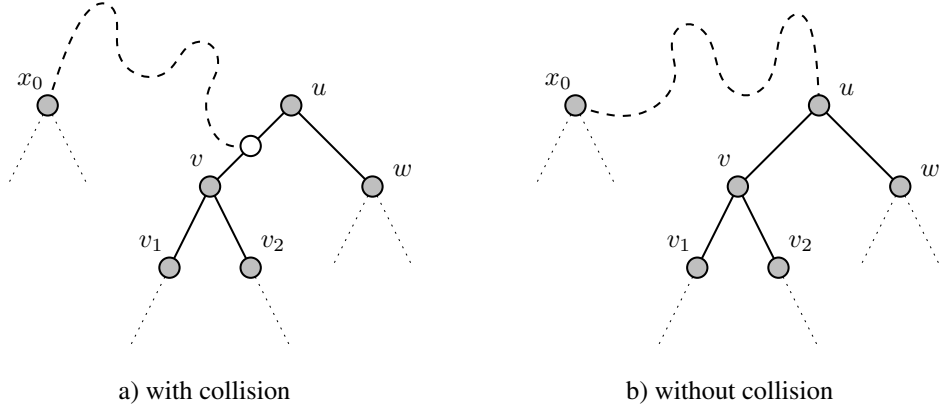\caption{Basic setup of the \iscollision\ routine in Figure~\ref{routine:iscollision}.
The subtrees $T_{x_0}$ and $T_u$ are edge-disjoint. The configuration
on the left shows a collision,
whereas the one on the right does not.
Note that in the collision case, the internal path of the
quartet $\{x_0,w,v_1,v_2\}$ is shorter than the edge $(u,v)$.
This is the basic observation behind the routine \iscollision.
We call $x_0$ the ``reference point''.}\label{fig:iscollision}
\end{center}
\end{figure}

An important step
of the reconstruction algorithm is to detect collisions, at least
when they are within a short distance.
Routine \iscollision\ defined in Figure~\ref{routine:iscollision}
and analyzed in Proposition~\ref{lem:iscollision} below performs
this task. We consider the following configuration, which we call
the ``basic collision setup.''
\begin{definition}[Basic Collision Setup]
We have two legal subtrees $T_{x_0}$ and
$T_{u}$ rooted at $x_0$ and $u$ in $T$.
We let $v, w$ be the children of $u$. We assume that
we are in either of the configurations depicted in Figure~\ref{fig:iscollision},
that is:
\begin{enumerate}
\item[a.] Either the path between $T_{x_0}$ and $T_u$
attaches in the middle of the edge $(u,v)$---in other words,
$T_{x_0}$ collides into $T_u$ at edge $(u,v)$;
\item[b.] Or it goes through
$u$---in other words, $T_{x_0}$ and $T_u$ are dangling.
\end{enumerate}
We call $x_0$ the ``reference point''.
\end{definition}
The purpose
of \iscollision\ is to distinguish between the two configurations above.
\begin{figure*}[h]
\framebox{
\begin{minipage}{16cm}
{\small \textbf{Algorithm} \iscollision\\
\textit{Input:} Four nodes $x_0$, $v$, $w$, $u$;
an edge length $h$;
a rooted forest and a distance matrix $(\forestno, \metricno)$;\\
\textit{Output:} \true\ or \false;
\begin{itemize}
\item \itemname{Children} Let $v_1, v_2$ be the children of $v$ in $\forestno$
(or $v_1 = v_2 = v$ if $v$ is a leaf);
\item \itemname{Internal\ Edge\ Length} Set
\begin{equation*} %\label{eq:defd}
\nu = \frac{1}{2} \left(
\metricno(v_1,x_0) +
\metricno(v_2,w) -
\metricno(v_1,v_2) -
\metricno(x_0,w)\right);
\end{equation*}
(set $\nu = +\infty$ if any of the distances is $+\infty$)
\item \itemname{Test} If $(h-\nu) > f/2$ return \true,
else return \false;
\end{itemize}
}
\end{minipage}
} \caption{Routine \iscollision.}
\label{routine:iscollision}
\end{figure*}
%\begin{figure}\centering
%\includegraphics[width=10cm, height=7cm]{FourPointGood.eps} \caption{Computing distances from
%closest common ancestor in {\sc IsSplit}.} \label{fig:fourpointgood}
%\end{figure}

\begin{proposition}[Correctness of \iscollision]\label{lem:iscollision}
Consider the Basic Collision Setup.
In particular, assume that one of the two configurations
in Figure~\ref{fig:iscollision} holds.
Let
$\quartet=\{x_0, v_1, v_2, w\}$ and let $\metricno$ and $h$ be the
parameters passed to \iscollision\ in Figure~\ref{routine:iscollision}. Assume that all edges in the tree satisfy
$d(e) \geq f$.
If
\begin{equation*}
\left|d(x,y) - \metricno(x,y)\right| < \frac{f}{8},\qquad \forall x,y \in \quartet,
\end{equation*}
and
\begin{equation*}
\left|h - d(u,v)\right| < \frac{f}{4},
\end{equation*}
then the routine \iscollision\ returns \true\ if and only if
$T_{x_0}$ collides into $T_u$ at edge $(u,v)$.
\end{proposition}
\begin{proof}
Whether or not there is a collision, the quantity $\nu$ computed
in the routine is the length of the internal path of the quartet
split $x_0w|v_1v_2$. If there is not a collision, then this path
is actually equal to the path corresponding to the edge $(u,v)$. 
Otherwise, the length
of the path is shorter than $d(u,v)$ by at least $f$ by assumption.
(See Figure~\ref{fig:iscollision}.)
The proof follows.
\end{proof}

%\subsection{Tree Operations}\label{sec:treeoperations}

%We collect in this subsection a few extra
%combinatorial subroutines.

\section{Reconstruction Algorithm} \label{sec:algvars}

%Below we proceed with an informal summary of the algorithm. The
%analysis proving its corretness and therefore Theorem~\ref{thm}
%is given in the next section.

We proceed with a formal description of the reconstruction algorithm.
A detailed example can be found in Section~\ref{sec:example}. The
reader may want to take a look at the example before
reading the details of the algorithm.

\subsection{Description of the Algorithm}

\paragraph{Cherry Picking.} Recall that in a $3$-regular tree a {\em cherry} is a pair of leaves
at graph distance $2$. Roughly speaking, our reconstruction
algorithm proceeds from a simple idea: it builds the tree one
layer of cherries at a time. To see how this would work, imagine that
we had access to a ``cherry oracle,'' that is, a function $C(u,v,T)$
that returns the parent of the pair of leaves $\{u,v\}$ if the
latter forms a cherry in the tree $T$ (and say $0$ otherwise).
Then, we could perform the following ``cherry picking'' algorithm:

{\small\tt
\begin{itemize}
\item Currently undiscovered tree: $T' := T$; \item Repeat until
$T'$ is empty,
\begin{itemize}
\item For all $(u,v) \in \lcal(T') \times \lcal(T')$, if\\ $w :=
C(u,v,T') \neq 0$, set $\parent{}(u):= \parent{}(v) := w$; \item
Remove from $T'$ all cherries discovered at\\ this step.
\end{itemize}
\end{itemize}
}
Unfortunately, the cherry oracle cannot be simulated from short
sequences at the leaves. Indeed, as we discussed in Section~\ref{sec:stattech},
short sequences provide only
``local'' metric information on the structure of the tree.
See the example in Section~\ref{sec:example} for an illustration
of the problems that may arise.
Nevertheless, the above scheme can be roughly followed by making a
number of modifications which we now describe.
% The full algorithm can be found in the Appendix.
% EM: Removed
%The description of the algorithm uses the following notation and
%conventions:
%\begin{itemize}

%\item A \emph{$g$-cherry} is a cherry where both edges have length less or equal to $g$.

%\item Let $M>0$. Let $T$ be a tree and $F$ be the subforest of $T$ where we keep all the
%leaves and only those nodes with the following property: they are on a
%path of length at most $M$ between two leaves of $T$.
%We say that a pair of leaves $\{u,v\}$ is an \emph{$M$-local
%  $g$-cherry} in $T$ if $\{u,v\}$ is a $g$-cherry in $F$ and there are
%at least two other leaves $u',v'$ s.t.
%$$\max\{d(u,u'),d(u,v'),d(v,u'),d(v,v')\} \leq M$$
%(the leaves $u',v'$ will act as ``witnesses'' of the cherry
%$\{u,v\}$).
% EM: Removed mentioning the appendix

%\item A \emph{pseudoleaf} is a current active node. The set of
%active nodes at iteration $i$ of the algorithm will be denoted by
%$\widehat{L}_i$ and will correspond to the leaves of the currently
%undiscovered part of the tree.

%\end{itemize}

%\paragraph{Blindfolded Cherry Picking.}
The high-level idea of the algorithm, which we call
\textsc{Blindfolded Cherry Picking} (BCP), is to apply the cherry
picking scheme above with two important differences:
\begin{itemize}
\item \itemname{Reconstructed\ Sequences}
Leaf sequences provide only local metric information ``around the leaves.'' To gain
information about higher, internal nodes of the tree, we reconstruct sequences at
the internal nodes of our partially reconstructed subforest and compute local metric
information ``around these nodes.'' By repeating this process, we gain information
about higher and higher parts of the tree.

\item \itemname{Fake\ Cherries} Moreover,
because of the local nature of our information, some of the
cherries we pick may in fact turn out \emph{not} to be cherries---that is,
they correspond to a path made of more than two edges in the true tree.
(See Section~\ref{sec:example} for an example.) As it turns out, this
only becomes apparent once a larger fraction of the tree is
reconstructed, at which point a subroutine detects the ``fake''
cherries and removes them.
\end{itemize}
The full algorithm is
detailed in Figures~\ref{routine:bcp},~\ref{routine:cherryid},~\ref{routine:fakecherry},
and~\ref{routine:bubble}.
We now explain its main components.
The parameters $k$, $\erada$, and $\tol$ will be set
in Section~\ref{sec:analysis}.

\begin{figure}
\begin{center}
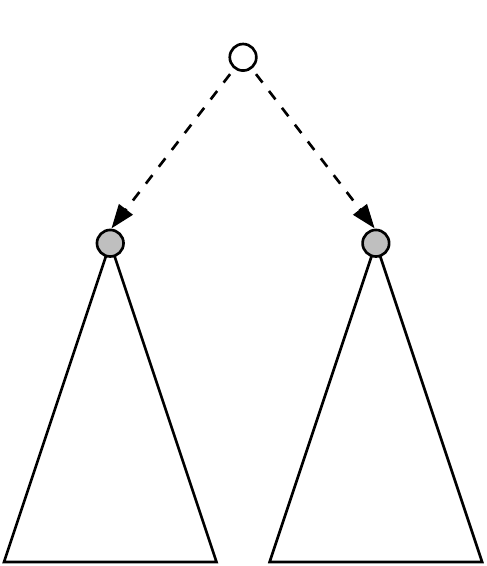\caption{Adding a cherry $(x,z,y)$ to $\forestno \supseteq
\{T_x, T_y\}$.}\label{fig:addingcherry}
\end{center}
\end{figure}

\paragraph{Adding a Cherry.} The algorithm BCP maintains
a partially reconstructed subforest of the true tree,
or more precisely, a legal subforest $\forestno$ of $T$.
The main operation we employ
to ``grow'' our partially reconstructed subforest is the
merging of two subtrees of $\forestno$ at their roots---an operation we call
``adding a cherry,'' in reference to the cherry picking
algorithm above.
%A {\em rooted cherry} is a pair
%of nodes $x,y$ with the same parent $z$.
Suppose the current
forest $\forestno$ contains two
edge-disjoint legal subtrees $T_x$ and $T_y$. We merge them
by creating a new node $z$ and adding the edges
$(z,x)$ and $(z,y)$ as in Figure~\ref{fig:addingcherry}.
We sometimes denote the pair of edges $\{(z,x),(z,y)\}$ by $(x,z,y)$.
We call this operation {\em adding cherry $(x,z,y)$
to $\forestno$}.
%\footnote{
%Typically the new edges form a cherry in the remaining
%unexplored tree.
%}

\paragraph{Identifying ``Local'' Cherries.}
As we explained above, we cannot hope to identify with certainty the
cherries of the unexplored part of the tree
from {\em short} sequences at the {\em leaves}.
Instead, we settle for detecting what we refer to as ``local'' cherries,
roughly, cherries in a ``local'' neighbourhood 
around the roots of the current reconstructed subforest.
More precisely, a ``local'' cherry is a pair of roots
of the current subforest that passes a series of tests
as detailed below.

To determine whether two roots $v_1, w_1$ of the current forest form
a ``local'' cherry,
our routine~\cherryid\footnote{In~\cite{DaMoRo:06}, the routine was called
\textsc{CherryID}.} in Figure~\ref{routine:cherryid} performs three tests:
\begin{enumerate}
\item \itemname{Short\ Distance} The nodes $v_1, w_1$
are at a short distance (roughly $2g$);
\item \itemname{Local\ Cherry} For all pairs of
roots $v_2, w_2$ at short distance (roughly $5g$),
the quartet
\begin{equation*}
\quartet = \{v_1, w_1, v_2, w_2\},
\end{equation*}
admits the split $v_1 w_1 | v_2 w_2$;
\item \itemname{Edge\ Lengths} The edges
connecting $v_1$, $w_1$ to their (hypothetical) parent
are short (roughly $g$).
\end{enumerate}
The routine \cherryid\ has three key properties, proved in Section~\ref{sec:analysis}:
\begin{enumerate}
\item \itemname{Edge\ Disjointness} It
preserves the {\em edge-disjointness} of the current forest;

\item \itemname{Restricted\ Subforest} It builds a forest that is always a {\em legal restriction}
of the true tree;

\item \itemname{Short\ Edges} It guarantees that all edges of the restricted
forest are {\em short} (smaller than $g'$).
\end{enumerate}
These properties are crucial for the proper operation of the algorithm
(in particular for the correctness
of routines such as \distmet\ (Figure~\ref{fig:distmet}) and \isshort\ (Figure~\ref{routine:isshort}) as seen in Section~\ref{sec:stattech}).

Finally, another important property of~\cherryid\ is that it is guaranteed to detect
\emph{true} cherries---at least those that are ``locally witnessed.''
We now define this notion more precisely.
For a distance matrix $\metric{}$ and a set of nodes $\ncal$,
we denote
\begin{equation*}
\maxmetric{}(\ncal) = \max\left\{\metric{}(x,y) : \{x,y\} \subseteq \ncal\right\}.
\end{equation*}
\begin{definition}[Witnessed Cherry]\label{def:localcherry}
%A \emph{$g$-cherry} is a cherry where both edges have length $\leq g$.
Let $\forestno$ be a forest with path metric $\metric{}$.
We say that a pair of leaves $\{u,v\}$ is an \emph{$M$-witnessed
cherry} in $\forestno$ if $\{u,v\}$ is a cherry in $\forestno$ and there are
at least two other leaves $u',v'$ s.t.
\begin{equation*}
\maxmetric{}(\quartet) \leq M,
\end{equation*}
where $\quartet = \{u,v,u',v'\}$
(the leaves $u',v'$ will act as ``witnesses'' of the cherry
$\{u,v\}$).
\end{definition}

\paragraph{Detecting Collisions.}
The merging of subtrees through ``local'' cherries that are
{\em not} actual cherries eventually results in {\em collisions}
between subtrees of the current forest such as in Figure~\ref{fig:remote}.
Such configurations are undesirable since they do not allow to
complete the reconstruction by simple merging of the subtrees at their roots.
Therefore, we seek to detect these collisions
using the routine~\footnote{In~\cite{DaMoRo:06}, the routine was called
\textsc{FakeCherry}.}
in Figure~\ref{routine:fakecherry}.

After adding a new layer of ``local'' cherries,
\fakecherry\ checks whether new collisions can be found.
For this, we use routine \iscollision\ in Figure~\ref{routine:iscollision} from Section~\ref{sec:combtech}.
We actually perform {\em two} such \iscollision\ tests for each candidate edge
in the target subtree $T_{u_1}$ (see Figure~\ref{routine:fakecherry}). This is done precisely
for the same reasons that we did a ``multiple test'' in routine
\distmet\ (Figure~\ref{fig:distmet}) in Section~\ref{sec:stattech}: It guarantees that at least
one test is performed under the Basic Disjoint Setup (Dangling), which is
needed for \iscollision\ to be correct. See Section~\ref{sec:analysis} for details.
%\item
%\fakecherry. This routine identifies ``fake'' cherries.
%The routine is rather involved -- see Figure~\ref{fig:fakecherry} for the
%algorithm
%and Section~\ref{sec:proof2} for proof of its correctness.
%\paragraph{Scanning Nodes in a Subtree.}
Also, to make sure that we are always in either
of the configurations in Figure~\ref{fig:iscollision},
we scan the nodes of the target subtree $T_{u_1}$ in
``reverse breath-first search (BFS) order'':
\begin{enumerate}
\item Order the vertices in subtree $T_{u_1}$ according
to breath-first search $v_1,\ldots,v_s$;
\item Scan the list in reverse order $v_s,\ldots,v_1$.
\end{enumerate}
This way, we never find ourselves ``above'' a collision in $T_{u_1}$ (because
otherwise we would have encountered the collision first).

\begin{figure}
\begin{center}
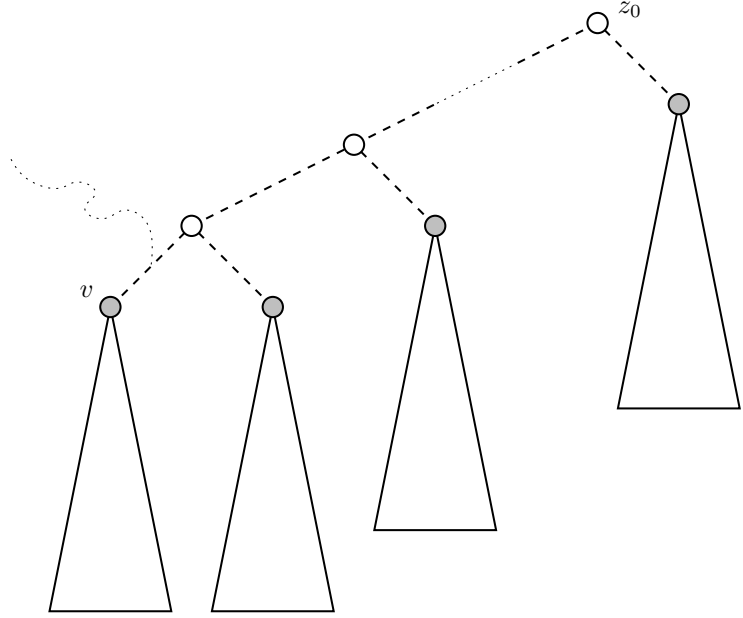\caption{The \bubble\ in Figure~\ref{routine:bubble} routine removes all
white nodes and dashed edges. The dotted curvy line indicates
the location of a collision to be removed.}\label{fig:bubble}
\end{center}
\end{figure}

\paragraph{Removing Collisions.}
Once a collision is identified we remove it using the
routine in Figure~\ref{routine:bubble}.\footnote{In~\cite{DaMoRo:06}, the routine was called \textsc{Bubble}.}
As seen in Figure~\ref{fig:bubble}, the routine essentially
removes the path from the collision all the way to the
root of the corresponding subtree.
For a
rooted forest $\forestno$, we use the notation
$\parent{\forestno}(x)$ to indicate the unique
parent of node $x$. 

We note that we actually go through the ``Collision Identification/Removal'' step {\em twice}
to make sure that we detect collisions occurring on edges adjacent
to the roots. See the proof of Lemma~\ref{lem:removal}.

We prove in Proposition~\ref{prop:progress} below
that, at each iteration of BCP, at least one cherry
is added that will not be removed at a later stage. Hence,
the algorithm makes progress at every iteration and eventually
recovers the full tree.

%Note that
%any non-root node $x$ has a unique ``sister'',
%i.e.~the other child of its parent. We denote
%the sister of $x$ by $\sister{\forestno}(x)$.

%Our proof of the main theorem is based on an explicit phylogenetic
%reconstruction algorithm from short sequences.
%The proof of its correctness can be found in the next section.
%A detailed example can be found in the next subsection. The
%reader may want to take a look at the example before
%reading the details of the algorithm.

%In Subsection~\ref{sec:example} we
%give an overview of the algorithm. Below we explain the role of
%the different variables/procedures.
\begin{figure*}[!ht]
\framebox{
\begin{minipage}{16cm}
{\small \textbf{Algorithm} \textsc{Blindfolded Cherry Picking}\\
\textit{Input:} Samples at the leaves $\{\sigma_u\}_{u\in [n]}$;\\
\textit{Output:} Estimated topology;
\begin{itemize}

\item {\bf 0) Initialization:}
\begin{itemize}
\item \itemname{Iteration\ Counter} $i:= 0$;
\item \itemname{Rooted\ Subforest} $\forest{0} := [n]$;
\item \itemname{Local\ Metric} For all $u, v \in \forest{0}$, set
$\metric{0}(u,v) = \distest(u,v; \forest{0}; \erada)$;
%\item \itemname{Roots\ of\ Current\ Subforest} $\roots{0} := [n]$;
%\item \itemname{Reconstructed\ Sequences} $\forall u\in [n],\, \recseq{u} := \sigma_{u}$;
\end{itemize}
%\item {\bf 1) Distance Estimation:}
%For all $(u,v) \in \widehat L_i \times \widehat L_i$, set
%$\hat d_i(u,v) := \mathrm{\textsc{DistEst}}(u,v)$ \textit{[see text]};

%\item {\bf 1) Sequence Reconstruction \& Distance Estimation:}
%\begin{itemize}
%\item For all $u, v \in \forest{i}$, set $\metric{i}(u,v) = \distest(u,v; \forest{i})$;
%\end{itemize}

\item {\bf 1) Local Cherry Identification:}
\begin{itemize}
\item Iteration: $i$;
\item Set $\forest{i+1} := \forest{i}$;
\item For all $(v_1,w_1) \in \binom{\roots{i}}{2}$,
\begin{itemize}
\item \itemname{Main\ Step} Compute $(\ischerry, l_v, l_w) := \cherryid\left( (v_1,w_1); (\forest{i}, \metric{i})\right)$;
\item If $\ischerry = \true$,
\begin{itemize}
\item \itemname{Update\ Forest} Create new node $u_1$ and add cherry
$(v_1,u_1,w_1)$ to $\forest{i+1}$;
%Add $N$ to $\roots{i+1}$ and $(v_1,N,w_1)$ to $\cherries{i+1}$;
%Remove $v_1, w_1$ from $\roots{i+1}$; Update forest $\forest{i+1}$
%accordingly, i.e.~add node $N$ and edges $(N, v_1), (N, w_1)$;
\item \itemname{Edge\ Lengths} Set $h(u_1,v_1) := l_v$ and $h(u_1, w_1) := l_w$;
%\item \itemname{Update\ Counter} Set $N:=N+1$;
\end{itemize}
\end{itemize}
\end{itemize}

%\item {\bf 2) Sequence Reconstruction:}
%\begin{itemize}
%\item For all $(u,w,v) \in \widehat C_{i}$, set $\hat\sigma_w :=
%\mathrm{\textsc{SeqRec}}(T^{\child}_{\leq w},w)$
%\textit{[see~(\ref{eq:seqrec})]};
%\end{itemize}

\item {\bf 2) Collision Removal:}
\begin{itemize}
\item \itemname{Update\ Metric} For all $x_1,x_2 \in \rho(\forest{i+1})$,
for all $u_\onetwo \in T^{\forest{i+1}}_{x_\onetwo}$, $\onetwo = 1,2$,
set
\begin{equation*}
\metric{i+1}(u_1,u_2) = \distmet(u_1,u_2; \forest{i+1}; \{h(e)\}_{e\in \forest{i+1}}; \erada; \tol);
\end{equation*}

\item Set $\forestno := \forest{i+1}$ and $\metric{} := \metric{i+1}$;
\item For all $(u_0,u_1) \in  \roots{}\times\roots{}$,
\begin{itemize}
\item Set $\hascollision:=\false$;
\item If $u_1$ is not a leaf,
\begin{itemize}
\item \itemname{Main\ Step} Compute $(\hascollision, z) := \mathrm{\fakecherry}\left((u_0,u_1); (\forest{}, \metric{})\right)$;
\end{itemize}
\item If $\hascollision = \mathrm{\true}$,
\begin{itemize}
\item \itemname{Update\ Forest} Compute $\forest{i+1} := \bubble(z; \forest{i+1})$;
\end{itemize}
\end{itemize}
\item \itemname{Second\ Pass} Set $\forestno := \forest{i+1}$ and repeat the previous step;
\end{itemize}

\item {\bf 3) Termination:}
\begin{itemize}
\item If $|\roots{i+1}| \leq 3$,
\begin{itemize}
\item Join
nodes in $\roots{i+1}$ (star if 3, single edge if 2);

\item Return (tree) $\forest{i+1}$;

\end{itemize}

\item Else, set $i:=i+1$, and go to Step 1.

\end{itemize}

\end{itemize}
}
\end{minipage}
} \caption{Algorithm \textsc{Blindfolded Cherry Picking}.}
\label{routine:bcp}
\end{figure*}

\begin{figure*}[!ht]
\framebox{
\begin{minipage}{16cm}
{\small \textbf{Algorithm} \cherryid\\
\textit{Input:} Two nodes $(v_1,w_1)$; current forest and distance matrix
$(\forestno, \metricno)$;\\
\textit{Output:} Boolean value and length estimates;
\begin{itemize}

\item Set $\ischerry:=\mathrm{\true}$ and $l_v=l_w=0$;

\item \itemname{Short\ Distance}
\begin{itemize}
\item If $\metricno(v_1,w_1)
> 2g + \tol$, then $\ischerry:=\mathrm{\false}$;
\end{itemize}

\item \itemname{Local\ Cherry}
%$(u_1,v_1) \in \widehat
%L_i\times \widehat L_i$ such that
%$u_1 < v_1$, $\{u_0,v_0\}\cap\{u_1,v_1\} =
%\emptyset$, and
%$\max\left\{\hat d_i(x_0,x_1)\,:\, x_\iota \in \{u_\iota, v_\iota\} \right\} \leq 5g+ \eps_2$.
%Then:
\begin{itemize}
\item Set $\neighbors = \left\{
(v_2, w_2) \in \binom{\rootsno}{2} :
%\max \left[\metricno(x,y) : x,y
%\in \{v_1,w_1,v_2,w_2\}\right]
%\left(\hat{d}(u_0,u_1),\hat{d}(u_0,v_1),\hat{d}(v_0,u_1),\hat{d}(v_0,v_1) \right)
\maxmetric{}(\{v_1,w_1,v_2,w_2\}) \leq 5g + \tol \right\} $;
\item If $\neighbors$ is empty, then
$\mathrm{IsCherry}:=\mathrm{\false}$;
Else, for all $(v_2,w_2) \in \neighbors$,
\begin{itemize}
\item If
$\mathrm{\issplit}\left((v_1,w_1),(v_2,w_2);\metricno\right)=\mathrm{\false}$
then set $\ischerry := \mathrm{\false}$ and \textit{break};
\end{itemize}
\end{itemize}

\item \itemname{Edge\ Lengths}
\begin{itemize}
\item If $\ischerry=\mathrm{\true}$,
\begin{itemize}
\item Let $x_1, x_2$ be the children of $v_1$ in $\forestno$ (or let $x_1 = x_2 = v_1$ if $v_1$ is a leaf);
\item Let $z_0$ be the closest node to $v_1$ in $\rootsno - \{v_1,w_1\}$ under $\metricno$;
\item Set $(b_v, l_v) := \mathrm{\isshort}\left((x_1,x_2), (w_1, z_0); \forestno; \erada; \tol/16 \right)$;
\item Repeat previous steps switching the roles of $v_1$ and $w_1$;
\item Set $\ischerry := b_v \land b_w$;
\end{itemize}
\end{itemize}

\item Return $(\ischerry, l_v, l_w)$;

\end{itemize}
}
\end{minipage}
} \caption{Routine \cherryid.} \label{routine:cherryid}
\end{figure*}

\begin{figure*}[!ht]
\framebox{
\begin{minipage}{16cm}
{\small \textbf{Algorithm} \fakecherry\\
\textit{Input:} Two roots $u_0,u_1$; directed forest and
distance matrix $(\forestno, \metricno)$;\\
\textit{Output:} Boolean value and node;
\begin{itemize}
%\item For $i =0,1$, set $T_i := T^{\child}_{\le u_i}$ and denote
%$C_i$ the set of cherries in $T_i$;
%\item
%Compute all pairwise distances $\hat d$ between $T_0$ and $T_1$
%using \textsc{DistEst} (some of these distances are actually
%wrong);
%\item $\forall (\kappa_0,\kappa_1) \in C_0\times C_1$ with
%$\kappa_i = (x_i,z_i,y_i)$, set $\hat d_M(\kappa_0,\kappa_1) =
%\max\{\hat d(v_0,v_1)\,:\,v_i\in\{x_i,y_i\}\}$; \item If $u_0$ is
%not in a cherry or $u_1$ is a leaf, \textsc{Break}.
%= 0,1$, unless $u_i$ is not in a cherry in $\widehat C_i$ or
%$u_{1-i}$ is a leaf, do
%\begin{itemize}
%\item Set $\mathrm{Stop}:=\mathrm{\textsc{false}}$; Let $\kappa_0
%= (x_0,u_0,y_0)$ be the cherry including $u_0$; Set $C' :=
%\{\kappa'\in C_1\,:\, \hat d_M(\kappa',\kappa) \leq 25g\}$;
%\textsc{Break} if empty; \item While $C' \neq \emptyset$ and
%$\mathrm{Stop}=\mathrm{\textsc{false}}$,
\item Set $\hascollision := \mathrm{\false}$ and $z := 0$;
\item Let $x_0, y_0$ be the children of $u_0$ in $\forestno$;
\item Scan through all nodes $v$ in $T_{u_1}$ (except $u_1$) in reverse
BFS manner,
\begin{itemize}
%\item Let $\kappa = (x, z, y) := \lowest(C')$; \item
\item Let $w := \sister{\forestno}(v)$ and $u:= \parent{\forestno}(v)$;
\item \itemname{Collision\ Test} Compute
\begin{equation*}
b_x := \mathrm{\iscollision}(x_0,v,w,u;h(u,v);(\forestno, \metricno)),
\end{equation*}
and
\begin{equation*}
b_y := \mathrm{\iscollision}(y_0, v, w, u; h(u,v); (\forestno, \metricno));
\end{equation*}
\item If $\hascollision := b_x \land b_y = \mathrm{\true}$
then set $z := v$ and \textit{break};
%the intersection of the triplet $\{x_0, x, y\}$; use the four
%point method on $\{x_r, x, y\}$ to compute the distance between
%$x$ and $w$, say $h$ (using a scheme similar to that in routine
%\textsc{DistEst}), check whether $h \neq \hat\gamma(x,z)$ (up to
%$2\eps_2$);
%\item \emph{[Collision Test 2]} Perform the previous
%step again with $y_r$ rather than $x_r$; \item If in both cases $h
%\neq \hat\gamma(x,z)$, then set
%$\mathrm{Stop}:=\mathrm{\textsc{true}}$ and set
%$w := z$; otherwise remove $\kappa$ from $C'$.
\end{itemize}
%\end{itemize}
%\item If $\mathrm{Stop}=\mathrm{\textsc{true}}$, perform
%$\mathrm{\textsc{Bubble}}(w, u_1)$.
\item Return $(\hascollision, z)$;
\end{itemize}
}
\end{minipage}
} \caption{Routine \fakecherry.} \label{routine:fakecherry}
\end{figure*}

\begin{figure*}[!ht]
\framebox{
\begin{minipage}{16cm}
{\small \textbf{Algorithm} \bubble\\
 \textit{Input:} Node $v$; rooted forest $\forestno$;\\
\textit{Output:} Rooted forest;
\begin{itemize}
%\item \itemname{Descendance} Let $\parent, \sister$
%be the descendance relationships implicit in $\fcal$;

\item If $v$ is not in $\forestno$
or $v$ is a root in $\forestno$, return $\forestno$;

\item Let $z_0$ be the root of the subtree of $\forestno$ in which $v$ lies;

\item Set $x:= v$;

\item While $x\neq z_0$,
\begin{itemize}
\item Set $x := \parent{\forestno}(x)$;
\item Remove node $x$ and its adjacent edges below it;
\end{itemize}
\item Return the updated forest $\forestno$.
\end{itemize}
}
\end{minipage}
} \caption{Routine \bubble.} \label{routine:bubble}
\end{figure*}

\paragraph{Implementation.}
We argue that the running time of the algorithm is $O(n^5 k)$ and that with the use of appropriate 
data structures it can be reduced further to $O(n^3 k)$. Let us start with the naive analysis. 
The distance estimations in the initialization step of the algorithm take overall time $O(n^2 k)$, 
since there are $O(n^2)$ pairs of leaves and each distance estimation between leaves takes linear 
time in the sequence length~$k$. Now, in every iteration of the algorithm:
\begin{itemize}
\item The Local Cherry Identification step takes overall time $O(n^4 k)$, since we consider $O(n^2)$ 
pairs of roots in the for-loop of this step, and each call of \cherryid\ requires $O(n^2 + nk)$ 
time---$O(n^2)$ time for the [Local Cherry] step and $O(n k)$ time for the [Edge Lengths] step, 
in which the \isshort\ call involves $O(1)$ ancestral sequence reconstructions on trees of size $O(n)$. 
\item The Collision Removal step requires $O(n^3 k)$ time in each iteration. Indeed, it performs $O(n^2)$ 
distance computations and each of these takes $O(n k)$ time, since it requires $O(1)$ ancestral 
sequence reconstructions on trees of size $O(n)$. It also performs $O(n^2)$ calls to \fakecherry\ 
and \bubble, each of which is a linear time operation on trees of size $O(n)$.
\end{itemize}
Since each iteration requires $O(n^4 k)$ time and there are $O(n)$ iterations (see proof in 
Section~\ref{sec:union}), the overall running time is $O(n^5 k)$. The above analysis is wasteful 
in 1) not reusing the already reconstructed ancestral sequences and in 2) performing various 
tests on pairs of nodes that are far apart in the tree. With the use of appropriate data structures, 
we could keep track of the ``local'' neighborhood of each node and restrict many computations 
to these neighborhoods. This should bring down the running time to $O(n^3 k)$ with a constant 
that would depend explicitly on $f, g$. The details are left to the reader.

\subsection{Example}\label{sec:example}

\begin{figure}
\begin{center}
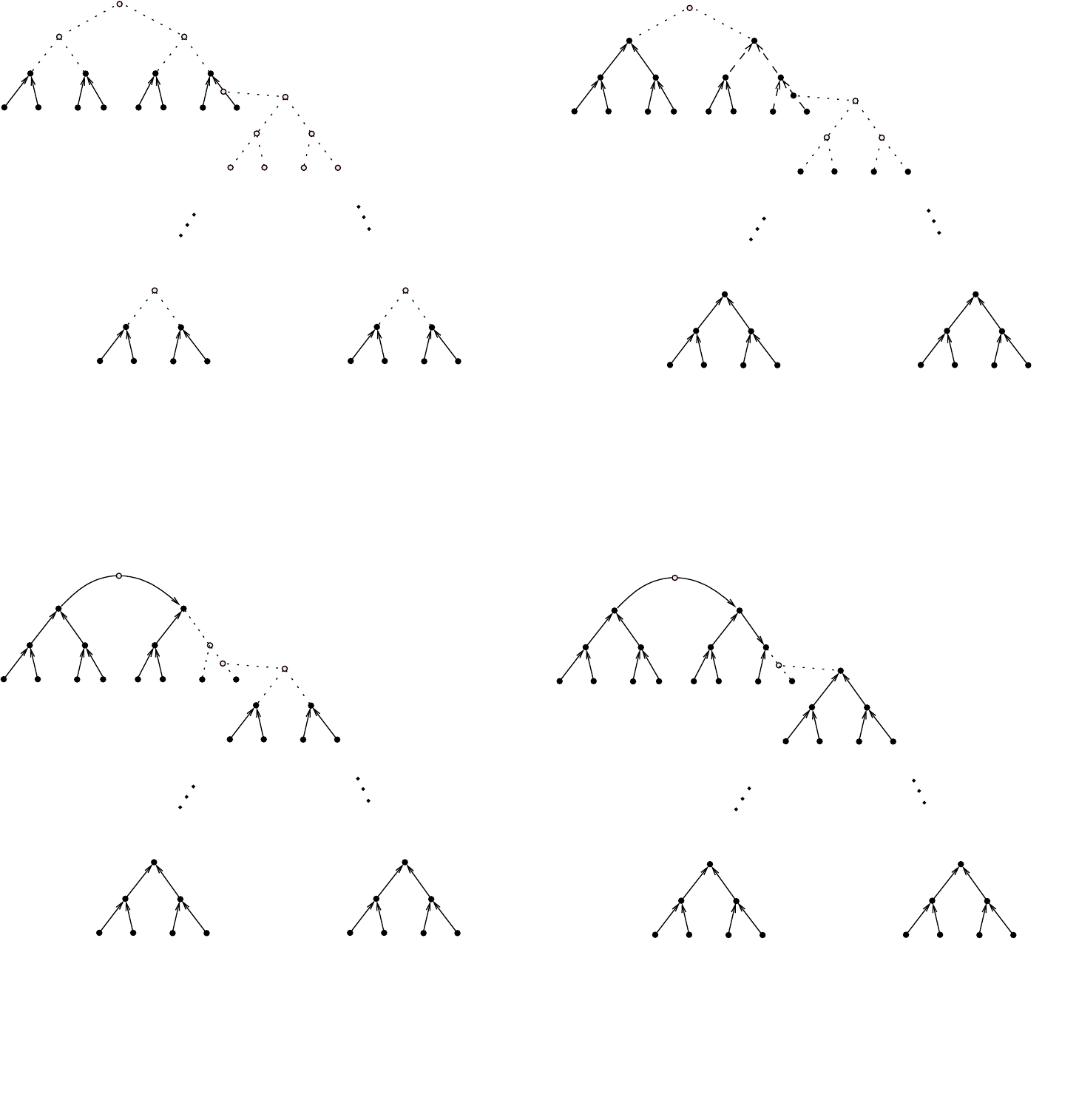\caption{Illustration of BCP's unraveling.}
\label{fig:full}
\end{center}
\end{figure}

We give a detailed example of the execution of the algorithm.
Consider the tree depicted in Figure~\ref{fig:full}a.
It is made of a large complete binary tree (on the right) with
a small $3$-level complete binary tree attached to its root (on the left).
All edges have length $g$, except $(v,x)$, $(x,v')$ and the edges adjacent to
the root of the small tree which have
length $g/2$.
%\begin{figure*}
%\centering
%\input{full.pstex_t}\caption{Illustration of BCP's unraveling.}
%\label{fig:full}
%\end{figure*}
In the figure, the subtree currently discovered by BCP is made of
solid arrows and full circles. The remaining (undiscovered) tree
is in dotted lines and empty circles. Assume that the length of
the sequences at the leaves allows us to estimate accurately
distances up to $5g$ (the actual constants used by the algorithm
can be found later).
%\begin{enumerate}
%\item[a)] \itemname{First\ Level}
\paragraph{First Level.} We first join into cherries all pairs of leaves that ``look''
like $g$-cherries in this ``local'' metric.
%i.e. they form a $g$-cherry
%5in all quartets of leaves. In particular, these selected pairs of leaves form
%$5g$-local $g$-cherries.
We are guaranteed to find all true
$g$-cherries.
However, consider pairs of leaves such as $u,v$
for which there is no ``local evidence''
that it does {\em not} form a cherry. Even though $u,v$ is not a cherry, it is joined
into a cherry by BCP.
Figure~\ref{fig:full}a depicts
the current forest after the first iteration of BCP.
%The top of the cherries are called pseudoleaves.
Before proceeding further, we apply our estimator $\ancname$ from Section~\ref{sec:ancestral}
to obtain reconstructed sequences at the roots of the current forest
and recompute the ``local'' metric.

%\item[b)] \itemname{Removing\ a\ Fake\ Cherry}
\paragraph{Removing a Fake Cherry.} We
subsequently proceed to join ``local'' $g$-cherries one
layer at a time, reconstructing internal sequences as we do so.
After many iterations, we find ourselves in the situation of
Figure~\ref{fig:full}b where most of the large complete tree has
been reconstructed (assume for now that edges $(u',v'')$,
$(v',v'')$, $(u,v')$, $(v,v')$ represented in dashed lines are
present). Now, the new information coming from sequences at
$y_1,\ldots,y_4$ provides evidence that $(u,v',v)$ is {\em not} a cherry
and that there is in fact a node $x$ on edge $(v,v')$. For
example, the quartet $\{y_1,y_2,u,v\}$ suggests that $u,v$ forms a
cherry with a $3g/2$-edge, which cannot hold in a true cherry. At
this point, we remove the ``fake'' cherry $(u,v',v)$ as well as
all cherries built upon it, here only $(u',v'',v')$. Note that we
have removed parts of the tree that were in fact reconstructed
correctly (e.g., the path between $u$ and $u'$).
% EM: Removed appendix

%\item[c)] \itemname{Rediscovering\ Removed\ Parts}
\paragraph{Rediscovering Removed Parts.} Subsequently,
BCP continues to join ``local'' cherries and
``rediscovers'' the parts of the tree that were removed. For instance,
in Figure~\ref{fig:full}c, the edge $(u',v'')$ is reconstructed again
but this time it forms a cherry with $(u'',v'')$ rather than $(v',v'')$.

%\item[d)] \itemname{Final\ Step}
\paragraph{Final Step.} Eventually, the full tree is correctly reconstructed
except may\-be for a few (at most $3$) remaining edges. Those can be added separately.
For example in Figure~\ref{fig:full}d only
the three edges around $x$ remain to be uncovered. Note that the
reconstructed tree has a root which is different from that of the
original tree.
% (i.e. the two trees differ as directed trees).
%\end{enumerate}

\section{Analysis I: Induction Step}\label{sec:analysis}

%{\red\bf [SR: Moved the progress stuff to the next section. The paragraph below is new.]}

In this section and the next, we establish that BCP reconstructs the phylogeny
correctly. In this section, we establish a number of combinatorial
properties of the current forest $\forest{i}$ grown by BCP.
Then, in the next section, we prove that the ``correctly reconstructed subforest'' of
$\forest{i}$ increases in size at every iteration.

\paragraph{Parameters.} Let $\delta > 0$ be the error tolerance.
Each application of Proposition~\ref{lem:statisticsRemote}
has an error of $O(n^{-\gamma})$. We will do $O(n^{3})$ distance
estimations so that by the union bound we require $O(n^{3-\gamma}) \leq \delta$.
Let $\tol < \frac{1}{8}\min\{f,g'-g\}$ and $\eradc > 6g$. Set $M > \eradc + 4g'$; 
take $\erada >  M+2\bias(g') + 4g'$; and choose $k$ to be equal to the maximum sequence length requirement for Proposition~\ref{lem:statisticsRemote} with parameters $\epsilon$, $M$ and $\erada$, for Proposition~\ref{lem:statisticsPair} with parameters $\epsilon$ and $M$, and for Proposition~\ref{lem:statisticsQuartet} with parameters $\tol/16$, $M$ and $\erada$.

%The following notation will be used in the proofs: $T$ is the
%phylogenetic tree that produced the data; $0 < f < g < +\infty$
%are lower and upper bounds on the length of every edge in $T$,
%where $g<g^*$; $k = c \log{n}$
% EM: s changed to k
%is the number of samples available at the leaves, where
%$c=c(f,g,\delta)$ is determined by the proof, $\delta$ being the
%probability of error. The proof yields a bound of $\exp(-c''(f,g)
%k) \leq n^{-\gamma}$ on the probability of failure for each
%distance estimation. We will have at most $n^{10}$ distance
%estimations and will therefore require that $n^{10-\gamma} \leq
%\delta$.

%Assume that $g$ satisfies the inequality $2e^{-2g}>1$, which
%defines the space of values of $g$ for which full reconstruction
%with $O(\log{n})$ samples at the leaves is not forbidden by
%\cite{Mossel:04a}. Also, fix the constant $\eps < f/2$ such that
%if $g'=g+\eps$ then $g'$ satisfies $2e^{-2g'}>1$.
%Fix $\eps_2<\eps/8$.

%\subsection{Proof of Main Theorem}\label{sec:proof1}

\paragraph{Induction Step.} The following proposition establishes a number of properties
of the forest grown by BCP.
We assume here that the conclusion of Proposition~\ref{lem:statisticsRemote}
holds at every iteration of the algorithm, until the tree is fully recovered.
In the next section, we will prove that the latter is indeed true with high probability.
%The proof can be found in Section~\ref{sec:proof2}.
%Recall
%that $\widehat{L}_i$ is the set of pseudoleaves (i.e. active nodes) at iteration
%$i$ (see the algorithm in Figure~\ref{fig:bcp} for a precise definition).
\begin{proposition}[Properties of $\forest{i}$]\label{prop:properties}
Denote by $\forest{i}=\left\{T_{
u}\,:\, u \in \roots{i}\right\}$ the current forest at the
beginning of \textsc{BCP}'s $i$-th iteration. Also, assume that the conclusion of Proposition~\ref{lem:statisticsRemote} holds for each call 
of the routine \distmet\ throughout the execution of the algorithm. 
Then, $\forall i\ge 0$ (until the algorithm stops),
\begin{enumerate}

\item\label{item:disjoint} $\mathrm{[Legal\ Subforest]}$ $\forest{i}$ is a legal
subforest of $T$;

\item\label{item:gedges} $\mathrm{[Edge\ Lengths]}$ $\forall u\in \roots{i}$,
$T_{u}$ has edge lengths at most $g'$;

\item\label{item:estimation} $\mathrm{[Weight\ Estimation]}$
The estimated lengths $\{h(e)\}_{e\in \ecal(\forest{i})}$ of the edges in
$\forest{i}$ are within $\tol/16$
of their right values;

\item\label{item:collisions} $\mathrm{[Collisions]}$ There is no collision within distance $\eradc$. (See Definition~\ref{def:collisions}.)

\end{enumerate}
\end{proposition}
\noindent\textbf{Proof of Proposition~\ref{prop:properties}:} The proof is by induction on $i$.\\

\noindent $i=0$: The set $\roots{0}$ consists of the leaves of $T$. The
claims are
therefore trivially true.\\

\noindent $i>0$: Assume the claims are true at the beginning of
the $j$-th iteration for all $j\leq i$. By doing a step-by-step analysis of the
$i$-th iteration, we show that the claims are still true at the
beginning of the $(i+1)$-st iteration.
%The following lemma follows
%from Proposition~\ref{lem:statisticsQuartetWithCollisions}.
%\begin{lemma}[Correctness of \textsc{DistEst}]\label{lem:smallDistancesCorrect}
%After the completion of step 1, for all $u,v \in \widehat{L}_i$:
%$$\hat{d}_i(u,v) \le 12g \vee d(u,v) \le 12g \Rightarrow |d(u,v)-\hat{d}_i(u,v)|<\eps_2.$$
%\end{lemma}
%\noindent\textbf{Proof:}
%From the induction hypothesis (Claim~\ref{item:collisions}), it follows that in the
%beginning of the $i$-th iteration there is no collision at
%distance $20g$. So the claim follows from
%Proposition~\ref{lem:statisticsQuartetWithCollisions}. (A small detail to note
%is that the sequences at the nodes of the forest were
%reconstructed in different steps of the algorithm. However, the
%subtrees that were used for the reconstruction of each node are
%exactly those in the statement of
%Proposition~\ref{lem:statisticsQuartetWithCollisions}.)
%$\blacksquare$

%{\bf\red [SR: Expanded definition of remaining forest below.]}
We first analyze the routine \cherryid\ (Figure~\ref{routine:cherryid}).
For a legal subforest $\forestno$ of $T$, we denote the ``remaining'' forest by
\begin{equation*}
\remforestno = T - \forestno.
\end{equation*}
%where we keep the nodes in $\roots{}$.
More precisely,
if $\forestno = \{T_1,\ldots,T_\alpha\}$
then $\remforestno$ is the
forest obtained from $T$ as follows:
\begin{enumerate}
\item Remove all edges of $T$ in the union of
the trees $T_1,\ldots,T_{\alpha}$.
In particular, for those edges of a $T_i$ representing
a path in $T$, we remove all corresponding edges of $T$.
\item The nodes of $\remforestno$ are all
the endpoints of the remaining edges of $T$.
All other nodes of $T$ are discarded.
\end{enumerate}
Note that the set $\remforestno$ is in
fact a subforest of $T$.
\begin{lemma}[Local Cherry Identification]\label{lem:fourPointIsCorrectFor4nodes}
Let $\forest{i}$ be the current forest at the
beginning of the $i$-th iteration.
Then we have the following:
\begin{itemize}

\item If $\{v_1,w_1\}$ is a $5g$-witnessed cherry
in $\remforest{i}$, then it passes all
tests in \cherryid;

\item If $\{v_1,w_1\}$ passes all the tests in \cherryid, then 
\begin{equation*}
d(u_1,v_1) \leq g+2\tol,
\end{equation*}
and
\begin{equation*}
d(u_1,w_1) \leq g+2\tol,
\end{equation*}
where $u_1$ is the parent of $\{v_1,w_1\}$
(as defined by \cherryid\ and the [Update Forest] step in the main loop).
Moreover the length estimates satisfy
\begin{equation*}
|d(u_1,v_1) - l_v| < \frac{\eps}{16},
\end{equation*}
and
\begin{equation*}
|d(u_1,w_1) - l_w| < \frac{\eps}{16}.
\end{equation*}

\item
If two pairs $\{v_1,w_1\}$ and $\{v_2,w_2\}$ both 
pass all tests in \cherryid, then it must be that the paths
$\path_T(v_1,w_1)$ and $\path_T(v_2,w_2)$ are non-intersecting.

\end{itemize}
\end{lemma}
\begin{proof}
Suppose first that $\{v_1,w_1\}$ is a $5g$-witnessed cherry
in $\remforest{i}$ with witness $\{v_2,w_2\}$. Since there is
no collision within $\eradc > 6g$ by assumption, we are in
the dangling case. Moreover, all edge weights below $\{v_1,w_1,v_2,w_2\}$ 
have been estimated
within $\tol/16$ and all corresponding edge weights (possibly
corresponding to paths) are $\leq g'$.
Therefore, by Propositions~\ref{lem:statisticsRemote}
and~\ref{lem:issplit}, \distmet\ (Figure~\ref{fig:distmet}) is accurate within $\tol$,
\issplit\ (Figure~\ref{routine:issplit}) returns the correct answer and
\isshort\ (Figure~\ref{routine:isshort}) is accurate within $\tol/16$.
Hence, $\{v_1,w_1\}$
passes the three tests in \cherryid.

Conversely, suppose $\{v_1,w_1\}$ passes all tests in \cherryid.
By our assumptions, when \distmet\ returns a finite value, it is
accurate within $\tol$. 
%Therefore we must have $d(v_1,w_1) \leq 2g + 2\tol$
%and in particular there is no collision between the subtrees
%rooted at $v_1$ and $w_1$.
%Also, there is at least one witness. For each such witness $\{v_2,w_2\}$
%we have
%\begin{equation*}
%\bar{d}(\{v_1,w_1,v_2,w_2\}) \leq 5g + 2\tol,
%\end{equation*}
%and, since
%there is no collision within $\eradc > 6g$, \issplit\ returns
%the correct answer. 
%This implies that $\{v_1,w_1\}$ is a cherry
%in $[\remforest{i}]_{5g+2\tol}$. 
Let $z_0$ be as in Figure~\ref{routine:cherryid}
(that is, $z_0$ is the closest node to $v_1$ in $\rootsno - \{v_1,w_1\}$ under 
the distorted metric).
The fact that $\{v_1,w_1\}$ previously passed the [Local Cherry] test
implies in particular that the diameter of $\{v_1,w_1,z_0\}$
must be less than distance $5g + 2\tol$.
In particular, there is no collision between
the subtrees rooted at $v_1$, $w_1$, and $z_0$ since $\eradc > 6g$.
Let $u_1$ be the intersection of $\{v_1,w_1,z_0\}$.
Then, \isshort\ returns an estimate within
$\tol/16 < \tol$ which in turn implies that $d(u_1,v_1)$ and $d(u_1,w_1)$ are $\leq g + 2\tol$.
%Hence, $\{v_1,w_1\}$ is a $(5g+2\tol)$-local $(g+2\tol)$-cherry.

For the third part, assume by contradiction that
the paths $\path_T(v_1,w_1)$ and $\path_T(v_2,w_2)$ intersect.
Then by the triangle inequality, the diameter of $\{v_1,w_1,v_2,w_2\}$
is at most $5g$. In particular, when \cherryid\ is applied to
$\{v_1,w_1\}$, the pair $\{v_2,w_2\}$ is considered in the
[Local Cherry] test and, since
there is no collision within $\eradc > 6g$, \issplit\ correctly returns
\false. That is a contradiction. 
\end{proof}

We can now prove Claims~\ref{item:disjoint},~\ref{item:gedges},
and~\ref{item:estimation} of Proposition~\ref{prop:properties}.
\begin{lemma}[Claims~\ref{item:disjoint},~\ref{item:gedges},
and~\ref{item:estimation}]
Let $\forest{i+1}$ be the current forest at the beginning of the $(i+1)$-st iteration. The Claims~\ref{item:disjoint},~\ref{item:gedges},
and~\ref{item:estimation} of the induction hypothesis hold for $\forest{i+1}$.
\end{lemma}
\begin{proof}
Since \bubble\ (Figure~\ref{routine:bubble}) only {\em removes} edges from the current forest,
it is enough to prove that after the completion of the {\em Local Cherry
Identification} step the
resulting forest satisfies Claims~\ref{item:disjoint},~\ref{item:gedges},
and~\ref{item:estimation}.

\noindent \textbf{Claim~\ref{item:disjoint}.}
By the induction hypothesis, $\forest{i}$ is legal. We only need to check that $\forest{i+1}$ is edge-disjoint.
Suppose on the contrary that the forest is not edge-disjoint. Also,
suppose that, along the execution of the {\em Local Cherry Identification} step,
the forest stops
being edge-disjoint when cherry $(v_1,u_1,w_1)$ is added to
$\forest{i}$. Then one of the following must be
true:

\begin{enumerate}
\item There is an ``old'' root $z \in \roots{i}$ such that $\path_T(v_1,w_1)$ is edge-sharing with
$T^{\forest{i}}_{z}$. But then there
is a collision in $\forest{i}$ within distance $2g+2\tol$ which contradicts the
induction hypothesis (Claim~\ref{item:collisions}).

\item There is a ``new'' root $z \in
\roots{i+1}\setminus \roots{i}$ with corresponding cherry $(x,z,y)$
such that $\path_T(v_1,w_1)$
is edge-sharing with
$\path_T(x,y)$.
In that
case $v_1w_1|xy$ is not the correct split
and
\begin{equation*}
\bar{d}(\{v_1,w_1,x,y\}) \leq 5g,
\end{equation*}
by the triangle inequality and Proposition~\ref{lem:statisticsRemote}.
But then, by Lemma~\ref{lem:fourPointIsCorrectFor4nodes},
\cherryid\ (Figure~\ref{routine:cherryid}) rejects $\{v_1,w_1\}$ when performing the [Local Cherry]
test with witness $\{x,y\}$---a contradiction.
\end{enumerate}

\noindent \textbf{Claim~\ref{item:gedges} and~\ref{item:estimation}.} Follows from the
induction hypothesis and Lemma~\ref{lem:fourPointIsCorrectFor4nodes}.
\end{proof}

It remains to prove Claim~\ref{item:collisions} of Proposition~\ref{prop:properties}. This follows from the following analysis of \fakecherry\ (Figure~\ref{routine:fakecherry}). Note that, since Claim~\ref{item:collisions} holds for all iterations $j\leq i$, it must be the case that any new collision between two trees involves at least one of the new edges of these trees added in the {\em Local Cherry Identification} step. We call an edge {\em deep} if it is not adjacent to a root in the current forest. Otherwise we call the edge a {\em top} edge. We show that the first pass of the {\em Collision Removal} step removes all collisions into deep edges. At the beginning of the second pass, all collisions (if any) must be into top edges. We show that the second pass cleans those up. We first prove that there are no false positives in the Collision Removal step.
\begin{lemma}[Collision Removal: No False Positive] \label{lem: no false positives}
Let $\forest{i+1}$ be the current forest at the beginning of the first or second pass of the Collision Removal step of the $i$-th iteration, and let $u_0, u_1 \in \roots{i+1,1}$ be the roots of the trees $T_0 = T^{\forest{i+1,1}}_{u_0}$ and $T_1 = T^{\forest{i+1,1}}_{
u_1}$. Let $v$ be a node in $T_1$, and suppose that $T_0$ does not collide into $T_1$ below $v$ or on the edge immediately above it. Then no collision is detected in the corresponding step of \fakecherry.
\end{lemma}
\begin{proof}
Let $x_0$, $y_0$ be the children of $u_0$. It suffices to show that either $b_x$ is \false\ or $b_y$ is \false\ (see Figure~\ref{routine:fakecherry}). Without loss of generality we can assume that the path connecting $u_0$ to $u_1$ does not pass through $x_0$. In particular, we are in the case b) of Figure~\ref{fig:iscollision}. If any of the distances passed to \iscollision\ is $+\infty$, \iscollision\ returns \false. Otherwise, by our assumption on the output of \distmet, the assumptions of Proposition~\ref{lem:iscollision} are satisfied. Hence, \iscollision\ returns \false\ in that case as well.
\end{proof}

\begin{lemma}[Collision Removal]\label{lem:removal}
The first and second passes of the Collision Removal step satisfy:
\begin{enumerate}
\item Let $\forest{i+1,1}$
be the current forest at the beginning
of the first pass of the Collision Removal step of the $i$-th iteration, and let $u_0,u_1 \in \roots{i+1,1}$.
Suppose $T_0 = T^{\forest{i+1,1}}_{u_0}$ collides into $T_1 = T^{\forest{i+1,1}}_{
u_1}$ within distance $\eradc$ on a deep edge $e=(u,v)$ of $T_1$. Then \fakecherry\ in Figure~\ref{routine:fakecherry} correctly detects the collision.

\item Let $\forest{i+1,2}$
be the current forest at the beginning
of the second pass of the Collision Removal step of the $i$-th iteration, and let $u_0,u_1 \in \roots{i+1,2}$.
Suppose $T_0 = T^{\forest{i+1,2}}_{u_0}$ collides into $T_1 = T^{\forest{i+1,2}}_{
u_1}$ within distance $\eradc$. Then \fakecherry\ correctly detects the collision.

\end{enumerate}
\end{lemma}
\begin{proof}
1. Denote by $x_\onetwo, y_\onetwo$ the children of $u_\onetwo$,
$\onetwo=0,1$.
Consider the Basic Collision Setup of Section~\ref{sec:combtech}.
By the remark above the statement of the lemma, the path
coming from $e$ enters $T_0$ through a top edge, or at $u_0$.
(See Figure~\ref{fig:fake}.)
\begin{figure}
\begin{center}
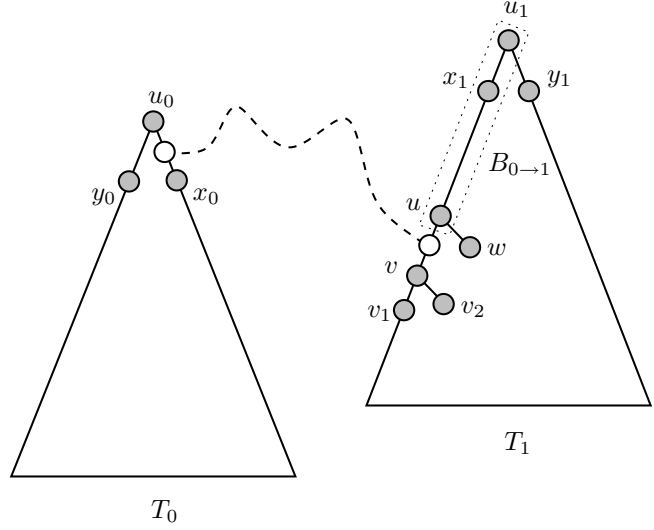\caption{Configuration in the proof of Lemma~\ref{lem:removal}.}\label{fig:fake}
\end{center}
\end{figure}
Let
\begin{equation*}
A_{0\to 1}=\left\{z\in \vcal(T_1):\text{$e$ is not in the subtree of $T_1$ rooted at $z$}\right\},
\end{equation*}
and
\begin{equation*}
B_{0 \to 1} = \vcal(T_1) \setminus A_{0 \to 1}.
\end{equation*}
Denote by $z$ the current-node variable used by \fakecherry\ as it scans the tree $T_1$ in a reverse BFS manner. Observe that, for all $z \in A_{0\rightarrow 1} - \{v\}$,
the Basic Collision Setup of Proposition~\ref{lem:iscollision} holds for both $x_0$ and $y_0$. Hence, by Lemma~\ref{lem: no false positives}, \iscollision\ in Figure~\ref{routine:iscollision} returns the correct answer. Furthermore, in the case of $v$, \iscollision\
returns \true\ for both $x_0$ and $y_0$ since the collision
is within $\eradc$ and
\begin{equation*}
\bar{d}(\{x_0,y_0,v_1,v_2,w\}) \leq \eradc + 4g' \leq M,
\end{equation*}
by Claim~\ref{item:gedges}. For $w$, \iscollision\ returns
\false\ because in that case $h(u,w) - \nu \leq - f + 3\eps \leq f/2$,
where $\nu$ is the estimated length of the internal edge of $\{x_0,w_1,w_2,v\}$
or $\{y_0,w_1,w_2,v\}$ with $w_1, w_2$ the children of $w$.
Finally, by scanning the nodes in reverse BFS order, we guarantee that
$v$ and $w$ are encountered before any node in $B_{0\to 1}$. Hence,
\fakecherry\ identifies correctly the collision on edge $e$.

2. Observe that, after a call to \bubble, the set of edges in the remaining forest is a subset of what it used to be. In particular, the subset of these edges involved in a collision decreases in size. Moreover, from the first part of the lemma, at the end of the first pass there is no collision remaining into deep edges. Given this, the argument above implies that any remaining collision within distance $\eradc$ will be found and removed in the second pass.
\end{proof}

\section{Analysis II: Tying It All Together}
\label{sec:union}

%{\red\bf [SR: This section is new but most of it is simply borrowed from the old Proposition 8.]}

In the previous section, we showed that provided Proposition~\ref{lem:statisticsRemote} holds at every iteration
the forest built is ``well-behaved.'' Below, we finish the proof of our main theorem by showing 
that Proposition~\ref{lem:statisticsRemote} indeed holds until termination and
that the algorithm eventually converges. We also discuss the issues involved
in extending $O(\log n)$-reconstruction beyond the $\quantum$-Branch Model in Section~\ref{section:beyond}.

\subsection{Quantifying Progress}

Our main tool to assess the progress of the algorithm is the notion
of a fixed subforest---in words, a subforest of the current forest
that will not be modified again under the normal operation of the algorithm.
%{\bf\red [SR: Clarified definition below.]}
\begin{definition}[Fixed Subforest] Let $\forestno$ be a legal
subforest of $T$. Let $u \in \vcal(\forestno)$. We say that $u$ is
\emph{fixed} if $T^{\forestno}_{u}$ is {\em fully reconstructed}, or in
other words, $T^{\forestno}_{u}$ can be obtained from $T$ by
removing (at most) one edge adjacent to $u$ and keeping the subtree containing $u$.
%{\blue EM: In the line above - it seems like if you remove the edge, the tree is determined uniquely - it consists all edges and vertices that are still connected to $u$}.
Note that
descendants of a fixed node are also fixed. We denote by
$\forestno^*$ the rooted subforest of $\forestno$ made of all fixed
nodes of $\forestno$. We say that $\forestno^*$ is the \emph{maximal fixed
subforest} of $\forestno$.
\end{definition}
Let
$\forest{i} = \{T_1,\ldots,T_\alpha\}$ be the current forest at iteration $i$ with remaining forest
%\begin{equation*}
$\remforest{i} = \{T'_1,\ldots,T'_\beta\}$.
Assuming the conclusion of Proposition~\ref{prop:properties} holds,
%\end{equation*}
% $\fcal'$ is the
%forest obtained from $T$ by removing all the edges in the union of
%the trees $T_1,\ldots,T_{\alpha}$. The nodes of $\fcal'$ are all
%the endpoints of the remaining edges. Since the trees
%$T_1,\ldots,T_{\alpha}$ are edge disjoint, the set $\fcal'$ is in
%fact a subforest of $T$.
each leaf $v$ in $\remforest{i}$ satisfies
exactly one of the following:
\begin{itemize}
\item \textbf{Fixed Root:} $v$ is a root of a fully reconstructed
tree $T_a \in \forest{i}$ (that is, $T_a$ is in
$\forest{i}^*$);
\item \textbf{Colliding Root:} $v$ is
a root of a tree $T_a \in \forest{i}$ that contains a collision
(that is, $T_a$ is not in $\forest{i}^*$);
\item \textbf{Collision Node:}
$v$ belongs to a path connecting
two vertices in $T_a \in \forest{i}$ but is
not the root of $T_a$ (that is, it lies in the ``middle'' of an edge of $T_a$).
\end{itemize}
Note in particular that the fixed roots of $\remforest{i}$ are roots in $\forest{i}^*$
(although not all roots of the maximal fixed subforest $\forest{i}^*$ 
are fixed roots in $\remforest{i}$ as they may lie below a colliding root).
We also need the notion of a fixed bundle---in words, a cherry (along with witnesses) that
will be picked by the algorithm at the next iteration and remain until termination.
\begin{definition}[Fixed Bundle]
A \emph{bundle in $\remforest{i}$} is a group of four leaves in $\remforest{i}$ such that:
\begin{itemize}
\item Any two leaves in the bundle are at topological distance
at most $5$ in $T$;
\item It includes at least one cherry of $\remforest{i}$.
\end{itemize}
A \emph{fixed
bundle} is a bundle in $\remforest{i}$ whose leaves are fixed roots.
\end{definition}
The following lemma is the key to our convergence argument.
It ensures that a fixed bundle always exists and hence that progress is made at every iteration.
\begin{proposition}[Existence of a Fixed Bundle]\label{prop:bundle}
Assume $\remforest{i}$ has at least $4$ leaves and satisfies the conclusion of Proposition~\ref{prop:properties}.
Then, $\remforest{i}$ contains at least one fixed bundle.
\end{proposition}
\noindent\textbf{Proof:}
%{\bf\red[SR: Expanded proof below.]}
To avoid confusion between the forests $\forest{i}$ and $\remforest{i}$, 
we will refer to the forest $\remforest{i}$ as the \emph{anti-forest},
to its trees as \emph{anti-trees} and to its leaves as \emph{anti-leaves}.
We first make a few observations:
\begin{enumerate}
\item The anti-forest $\remforest{i}$ is binary, that is, all its nodes have degrees in $\{1,3\}$.
Indeed, note first that
since $T$ is binary, one cannot obtain nodes of degree higher than $3$ by removing
edges from $T$. Assume by contradiction that there is a node, $u$, of degree $2$ in
$\remforest{i}$. Let $w_1,w_2,w_3$ be the neighbors of $u$ in $T$ and assume that
$w_3 \notin \remforest{i}$. By construction of the anti-forest $\remforest{i}$ the edge $(u,w_3)$ is in the forest $\forest{i}$
(possibly included in an edge of $\forest{i}$ corresponding to a path of $T$).
Moreover, $u$ is a node of degree $1$ in $\forest{i}$. But this is a contradiction because the only
nodes of degree $1$ in a legal subforest of $T$ are the leaves of $T$ and $u$ cannot be a leaf of $T$.

\item\label{obs:2} A binary tree (or anti-tree) $T_0$ with $4$ or more leaves (or anti-leaves)
contains at least one bundle.
Indeed, let $L_0$ be the leaves of $T_0$. We call $R_0 = L_0$ the \emph{level-$0$ leaves} of $T_0$.
Now remove all cherries from $T_0$ and let $R_1$ be the roots of the
removed cherries---the \emph{level-$1$ leaves}.
Denote by $T_1$ the tree so obtained and note that its leaves $L_1$ contain $R_1$ as well as
some remaining level-$0$ leaves. Consider the cherries
of $T_1$ (there is at least one). If any such cherry is made of two level-$1$ leaves, then
the corresponding (descendant) level-$0$ leaves of $T_0$ form a bundle in $T_0$ and we are done.
(In that case, the diameter of the bundle is $4$.)
Suppose there is no such cherry. Note that there is no cherry in $T_1$ formed from
two level-$0$ leaves as those were removed in constructing $T_1$. Hence,
all remaining cherries of $T_1$ must contain at least one level-$1$ leaf. Now, merge
all such cherries to obtain $T_2$. Denote by $R_2$ the roots of the removed cherries---the
\emph{level-$2$ leaves}. Once again, by the argument above, any cherry of $T_2$ contains at least
one level-$2$ leaf. Any such cherry (there is at least one) provides a bundle by considering its descendant
level-$0$ leaves in $T_0$. 
(If the second leaf in the cherry is level-$0$, the diameter of the bundle is $4$.
If it is level-$1$ or level-$2$, the diameter is $5$.)
This proves the claim.

\item By Claim~\ref{item:collisions} in
Proposition~\ref{prop:properties}, there is no collision within $\eradc > 6g$ 
(see Definition~\ref{def:collisions}).
In particular,
collision nodes are at distance
at least $\eradc > 6g$ from any other anti-leaf (collision nodes, fixed roots, colliding roots)
in $\remforest{i}$. Therefore, if an
anti-tree in $\remforest{i}$ contains a collision node, then it has $\geq 4$
anti-leaves and, from the previous observation, it contains at least one
bundle. Moreover, this bundle cannot contain a collision node
since in a bundle  all anti-leaves are at distance at most $5g$ and collision nodes
are at distance at least $\eradc > 6g$ from all other anti-leaves in $\remforest{i}$.

\item\label{obs:1} From the previous
observations, we get the following: if a tree in $\remforest{i}$ contains
a collision, then either it has a fixed bundle or it has at least
one colliding root.
\end{enumerate}

We now proceed with the proof. Assume first that there is no
collision node in $\remforest{i}$. Then, there cannot be any colliding root either
because by definition a colliding root is the root of a tree containing a collision.
In particular, $\remforest{i}$ is composed of a single anti-tree whose
anti-leaves are all fixed roots.
Then, since by assumption $\remforest{i}$ has at least $4$ anti-leaves,
there is a fixed bundle by Observation~\ref{obs:2} above.

Assume instead that there is a collision node. Let $\widetilde{T}^{0}$ be
an anti-tree in $\remforest{i}$ with such a collision node, say $c^0$.
Then by Observation~\ref{obs:1}, either 1) $\widetilde{T}^{0}$
has a fixed bundle in which case we are done, or 2) one of $\widetilde{T}^{0}$'s anti-leaves is a colliding root,
say $r^0$. In the latter case, let $T^0$ be the tree in $\forest{i}$
whose root is $r^0$. The tree $T^0$ contains at least one collision node (included in an edge corresponding to a path of $T$).
This collision node, say $c^1$, is also contained in an anti-tree in $\remforest{i}$, say $\widetilde{T}^{1}$. Repeat the
argument above on $\widetilde{T}^{1}$, and so on.

%{\blue EM: I am still a bit lost. What does "share" above means? Below I think it will be good if we can define a path via the entering and exiting process and show that this path is locally simple 
%(the new part added is simple and does not intersect the last part added. Since this is a tree - this will prove the process terminates}

We claim that this process generates a simple path $P$ in $T$ starting from the node $c^0$ above, 
passing through an alternating sequence of colliding roots and collision nodes $r^0, c^1, r^1, c^2,\ldots$, and eventually
reaching a fixed bundle.
%Indeed in each step where a fixed bundle is \emph{not} found we ``exit'' a tree in $\forest{i}$
%via a node that is not
%its root (the collision node $u_b$ above) and ``enter'' a new tree
%in $\forest{i}$ at its root (the colliding root $v_b$ above).
Indeed, since there is no cycle in $T$, the only way for $P$ not to be simple is for it
to ``reverse on itself.'' But note that by definition, for all $j$ we have $c^j \neq r^j$ and
$r^j \neq c^{j+1}$. Moreover the simple paths $c^j \to r^j$ and $r^j \to c^{j+1}$ belong
respectively to the anti-forest $\remforest{i}$ and the forest $\forest{i}$ (possibly with subpaths
collapsed into edges). In particular, their edges (in $T$) cannot intersect. Hence, $P$ is simple. 
Finally, since $T$ is finite, this
path cannot be infinite, and we eventually find a fixed
bundle.
$\blacksquare$

\subsection{Proof of the Main Theorem}

Consider now the $\quantum$-Branch Model. The proof of convergence works by considering
first the hypothetical case where all distance estimates computed by the algorithm are perfectly accurate,
%conclusion of 
%Proposition~\ref{lem:statisticsRemote} holds with $\eps = 0$, 
that is
the case where we have ``perfect'' local information. We denote by $\{\forestz{i}\}_{i\geq 0}$ the sequence
of forests built under this assumption. Note that, up to arbitrary choices (tie breakings, orderings, etc.), this sequence
is \emph{deterministic}. We now show that it terminates in a polynomial number of steps with the correct tree.
\begin{proposition}[Progress Under Perfect Local Information]\label{prop:progress}
Assume $\forestz{i}=\left\{T_{
u}\,:\, u \in \rootsz{i}\right\}$ is the current forest at the
beginning of \textsc{BCP}'s $i$-th iteration under perfect local information
with corresponding maximal fixed subforest $\forestz{i}^*$.
Then for all $i\geq 0$ (before the termination step),
$\forestz{i}^*\subseteq\forestz{i+1}^*$ and
$|\vcal(\forestz{i+1}^*)| > |\vcal(\forestz{i}^*)|$.
%In particular, the algorithm terminates after $O(n)$ iterations
%with the correct tree.
\end{proposition}
\noindent\textbf{Proof of Proposition~\ref{prop:progress}:} We
first argue that
$\forestz{i}^*\subseteq\forestz{i+1}^*$. Note that the
only routine that removes edges is \bubble\ in Figure~\ref{routine:bubble}. Since \bubble\ only removes edges above identified collisions and $\forestz{i}^*$ is fully
reconstructed, it suffices to show that collisions identified by
\fakecherry\ in Figure~\ref{routine:fakecherry} are actual collisions. This follows from Lemma~\ref{lem: no false positives}.

%or lie ``above'' an
%actual collision---i.e. are on a cherry located on the path
%between the actual collision and the root. Indeed, since
%\textsc{Bubble} removes only edges ``above'' presumed collisions,
%this would then imply that no edge in $\widehat\fcal_{i}^*$ can be
%removed. We now prove the claim by analyzing the behavior of
%\textsc{FakeCherry}. We use the notation defined in the routine.
%Consider the collision tests in \textsc{FakeCherry}. The key point
%is to observe the following:
%\begin{itemize}
%\item if cherry $\kappa=(x,z,y)$ is in $\widehat\fcal_{i}^* \cap
%T_{i}$ and $\kappa_{1-i}=(x_r,u_{1-i},y_r)$, then at
%least one of $x_r$ or $y_r$ has a reconstruction bias that is
%independent from the bias at both $x$ and $y$; therefore this
%``correct'' witness will not observe a collision (using
%Proposition~\ref{lem:statisticsQuartetWithCollisions}); \item if cherry $\kappa$ is in
%$(T-\widehat\fcal_{i}^*) \cap T_{i}$, then all the cherries
%above $\kappa$ (on the path to $u_{i}$) cannot be in
%$\widehat\fcal_{i}^*$ and therefore applying \textsc{Bubble} to
%$\kappa$ does not modify $\widehat\fcal_{i}^*$.
%\end{itemize}
%This proves the first claim.

To prove $|\vcal(\forestz{i+1}^*)| > |\vcal(\forestz{i}^*)|$, assume
$\forestz{i} = \{T_1,\ldots,T_\alpha\}$ with remaining forest
%\begin{equation*}
$\remforestz{i} = \{T'_1,\ldots,T'_\beta\}$.
From Proposition~\ref{prop:bundle}, it follows that $\remforestz{i}$ contains at least one fixed
bundle.
This immediately implies the second claim.
Indeed, by Lemma~\ref{lem:fourPointIsCorrectFor4nodes}
note that the cherry in the fixed bundle is found by
\cherryid\ in Figure~\ref{routine:cherryid} during the 
$(i+1)$-st iteration and is not removed by
the {\em Collision Removal} step from the argument above.
$\blacksquare$

\noindent\textbf{Proof of Theorem~\ref{thm}:} 
To prove our main theorem under the $\quantum$-Branch Model,
we modify our reconstruction algorithm slightly by rounding 
the estimates in Proposition~\ref{lem:statisticsRemote}
to the closest multiple of $\quantum$. Also, we choose a number of samples large enough so that
the distance estimation error is smaller than $\quantum/3$. 
In that case, \emph{we simply mimic the algorithm under perfect local information}.
Note that by Proposition~\ref{prop:progress} there are at most $O(n)$
iterations until termination under perfect local information. 
By a union bound, it follows that Proposition~\ref{lem:statisticsRemote}
holds true for all pairs of subtrees in $\{\forestz{i}\}_{i\geq 0}$
with high probability.

We can now conclude the proof.
By Proposition~\ref{prop:properties}, the current forest at each iteration is correctly
reconstructed. 
By Proposition~\ref{prop:progress} after $O(n)$
iterations there remain at most three nodes in $\roots{i}$
and at that point, from Proposition~\ref{prop:properties} Claim~\ref{item:collisions},
we have that $\forest{i}\equiv \forest{i}^{*}$. Therefore
the remaining task is to join the remaining roots and there
is only one possible topology. So when the BCP algorithm
terminates, it outputs the tree $T$ (as an undirected tree) with
high probability.
%and all estimated edges are within $\tol/16$ of
%their correct value.
%This concludes the proof. 

The tightness of
the value $g^{\ast} = \frac{\ln 2}{4}$ is justified by the
polynomial lower bound \cite{Mossel:04a} on the number of required
characters if the mutation probability $p$ on all edges of the
tree satisfies $2(1-2p)^2<1$. $\blacksquare$

\subsection{Beyond the $\quantum$-Branch Model?}\label{section:beyond}

%{\red\bf [SR: Added this.]}

Extending Theorem~\ref{thm} to \emph{continuous} edge lengths appears far from trivial.
The issue arises in the final union bound over all applications of Proposition~\ref{lem:statisticsRemote} 
(more specifically, the ancestral state reconstruction step) 
which is valid except with inverse polynomial probability over
the generated sequences (for a fixed subtree of $T$). Indeed, note that in general there are \emph{super-polynomially}
many restricted subtrees of the true tree $T$ where all edges (paths in $T$) are shorter than $g^\ast$.
Therefore, using only a simple union bound, we cannot hope to guarantee that ancestral state reconstruction is successful \emph{simultaneously}
on all possible partial reconstructed subtrees.

In the previous subsection, 
we avoided this problem by showing that under the $\quantum$-BM assumption the algorithm follows
a \emph{deterministic} realization path of polynomial length. Moving beyond this proof would require
proving that the ancestral state reconstruction can be performed on the \emph{random} forests generated by the algorithm---but
this is not straightforward as the partially reconstructed forest is generated by the same data that is used to perform the
ancestral estimation. We conjecture that the correlation created by this ``bootstrap'' is mild enough to allow our algorithm
to work in general, but we cannot provide a rigorous proof at this point.

We remark that Mossel's earlier result~\cite{Mossel:04a} on the balanced case with continuous edge lengths
is not affected by the issue above because, there,
the reconstruction of the tree occurs one level at a time (there is no collision). Hence, ancestral state reconstruction is performed
only on \emph{fully} reconstructed subtrees---of which there are only polynomially many.

Note finally that, even under the discretization assumption made in this paper, achieving $\log(n)$-re\-cons\-truction
is nontrivial and does not follow from previous techniques. In particular, it can be shown that all previous rigorous
reconstruction algorithms for general trees
require polynomial sequence lengths even when all edge lengths are identical and below $g^\ast$.
See~\cite{Roch:08} for a formal argument of this type.

\section{Conclusion}

The proof of Steel's Conjecture~\cite{Steel:01} provides tight results
for the phylogenetic reconstruction problem. However, many theoretical and
practical questions remain:
\begin{itemize}
\item Can the discretization assumption be removed? We conjecture that the answer is yes.
\item
Can the results be extended to other mutation models? 
In subsequent work, Roch~\cite{Roch:09}
showed that our results hold for so-called
General Time-Reversible (GTR) mutation models
below the Kesten-Stigum bound.
Can the results
be extended to deal with ``rates across sites''
(see e.g.~\cite{Felsenstein:04})?
\item
What is the optimal $g$-value for the Jukes-Cantor model?
Is it identical to the critical value $g_{q=4}$ of the reconstruction problem
for the so-called ``Potts model'' with $q=4$? We note that it is a
long standing  open problem to find $g_{q=4}$.
The best bounds known are given in~\cite{MosselPeres:03,MaSiWe:04b}.
See also~\cite{MezardMontanari:06, Sly:08}.
\end{itemize}

\section*{Acknowledgments}

%C.D. is supported by CIPRES (NSF ITR grant \# NSF EF 03-31494).
%E.M. is supported by a
%Miller fellowship in Statistics and Computer Science, by a Sloan
%fellowship in Mathematics and by NSF grants DMS-0504245 and DMS-0528488.
%S.R. is supported by CIPRES (NSF ITR grant \# NSF EF 03-31494), FQRNT,
%NSERC and a Lo\`eve Fellowship.
S.R. thanks Martin Nowak and the Program for Evolutionary Dynamics
at Harvard University where part of this work was done. C.D. and
E.M. thank Satish Rao for interesting discussions. E.M. thanks
M.~Steel for his enthusiastic encouragement for studying the connections
between the reconstruction problem and phylogenetics.
We thank Allan Sly for helpful discussions.

\bibliographystyle{alpha}
\bibliography{thesis}

\end{document}